\documentclass[twoside,11pt]{article}
\usepackage{jmlr2e}
\RequirePackage{amsthm,amsmath,amssymb,graphicx}
\RequirePackage[colorlinks,citecolor=blue,urlcolor=blue]{hyperref}
\usepackage{graphicx,epstopdf}
\usepackage{natbib}
\usepackage{amsmath,amsthm}
\usepackage{booktabs}
\usepackage{paralist}
\usepackage{dsfont}
\usepackage{pdfpages}
\usepackage{wrapfig}
\usepackage[mathscr]{euscript}
\usepackage{filecontents}

\def\edt{\end{document}}
{\catcode `\@=11 \global\let\AddToReset=\@addtoreset}
\AddToReset{equation}{section}

\newtheorem{Theorem}{Theorem}
\newtheorem{Lemma}{Lemma}
\AddToReset{Theorem}{section}
\AddToReset{Lemma}{section}

\numberwithin{equation}{section}
\numberwithin{table}{section}
\numberwithin{figure}{section}

\newcommand{\junk}[1]{{}}
\DeclareMathOperator*{\argmin}{arg\,min}
\DeclareMathOperator*{\argmax}{arg\,max}

\ShortHeadings{JPEN Estimation of Covariance and Inverse Covariance Matrix}{Ashwini Maurya}
\firstpageno{1}
\begin{document}
\title{A Well-Conditioned and Sparse Estimation of Covariance and Inverse Covariance Matrices Using a Joint Penalty}

\author{\name Ashwini Maurya \email mauryaas@msu.edu \\
       \addr Department of Statistics and Probability\\
       Michigan State University\\
       East Lansing, MI 48824, USA}

\editor{}

\maketitle

\begin{abstract}
We develop a method for estimating well-conditioned and sparse covariance and inverse covariance matrices from a sample of vectors drawn from a sub-gaussian distribution in high dimensional setting. The proposed estimators are obtained by minimizing the quadratic loss function and joint penalty of $\ell_1$ norm and variance of its eigenvalues. In contrast to some of the existing methods of covariance and inverse covariance matrix estimation, where often the interest is to estimate a sparse matrix, the proposed method is flexible in estimating both a sparse and well-conditioned covariance matrix simultaneously. The proposed estimators are optimal in the sense that they achieve the minimax rate of estimation in operator norm for the underlying class of covariance and inverse covariance matrices. We give a very fast algorithm for computation of these covariance and inverse covariance matrices which is easily scalable to large scale data analysis problems. The simulation study for varying sample sizes and variables shows that the proposed estimators performs better than several other estimators for various choices of structured covariance and inverse covariance  matrices. We also use our proposed estimator for tumor tissues classification using gene expression data and compare its performance with some other classification methods.
\end{abstract}

\section{Introduction}

	With the recent surge in data technology and storage capacity, today's statisticians often encounter data sets where sample size $n$~is small and number of variables $p$ is very large: often hundreds, thousands and even million or more. Examples include gene expression data and web search problems [\citet{Clark1:7}, \citet{pass:21}]. For many of the high dimensional data problems, the choice of classical statistical methods becomes inappropriate for making valid inference. The recent developments in asymptotic theory deal with increasing $p$ as long as both $p$ and $n$ tend to infinity at some rate depending upon the parameters of interest.
	
The estimation of covariance and inverse covariance matrix is a problem of primary interest in multivariate statistical analysis. Some of the applications include: \textbf{(i)} Principal component analysis (PCA) [\citet{Johnstone:14}, \citet{Zou:37}]:, where the goal is to project the data on ``best" $k$-dimensional subspace, and where best means the projected data explains as much of the variation in original data without increasing $k$. \textbf{(ii)} Discriminant analysis [\citet{Mardia:19}]:, where the goal is to classify observations into different classes. Here estimates of covariance and inverse covariance matrices play an important role as the classifier is often a function of these entities. \textbf{(iii)} Regression analysis: If interest  focuses on estimation of regression coefficients with correlated (or longitudinal) data, a sandwich estimator of the covariance matrix may be used to provide standard errors for the estimated coefficients that are robust in the sense that they remain consistent under mis-specification of the covariance structure. \textbf{(iv)} Gaussian graphical modeling [\citet{Mein:20}, \citet{Wainwright:27}, \citet{Yuan:34},\citet{Yuan1:35}]:, the relationship structure among nodes can be inferred from inverse covariance matrix. A zero entry in the inverse covariance matrix implies conditional independence between the corresponding nodes.  \\

The estimation of large dimensional covariance matrix based on few sample observations is a difficult problem, especially when $n \asymp p$ (here $a_n \asymp b_n$ means that there exist positive constants $c$ and $C$ such that $c \le a_n/b_n \le C $). In these situations, the sample covariance matrix becomes unstable which explodes the estimation error. It is well known that the eigenvalues of sample covariance matrix are over-dispersed which means that the eigen-spectrum of sample covariance matrix is not a good estimator of its population counterpart [\citet{Marcenko:18}, \citet{Karoui1:16}]. To illustrate this point, consider $\Sigma_p=I_p$, so all the eigenvalues are $1$. A result from [\citet{Geman:12}] shows that if entries of $X_i$'s are i.i.d (let $X_i$'s have mean zero and variance 1)  with a finite fourth moment and if $p/n \rightarrow \theta <1 $, then the largest sample eigenvalue $l_1$ satisfies:
\begin{eqnarray*}
l_1~ \rightarrow ~(1+\sqrt{\theta})^2, ~~~~~~ a.s
\end{eqnarray*}
This suggests that $l_1$ is not a consistent estimator of the largest eigenvalue $\sigma_1$ of population covariance matrix. In particular if $n=p$ then $l_1$ tends to $4$ whereas $\sigma_1$ is $1$. This is also evident in the eigenvalue plot in Figure 2.1. The distribution of $l_1$ also depends on the underlying structure of the true covariance matrix. From Figure 2.1, it is evident that the smaller sample eigenvalues tend to underestimate the true eigenvalues for large $p$ and  small $n$. For more discussion on this topic, see \citet{Karoui1:16}. 

To correct for this bias, a natural choice would be to shrink the sample eigenvalues towards some suitable constant to reduce the over-dispersion. For instance, \citet{Stein:28} proposed an estimator of the form $\tilde{\Sigma}=\tilde{U} \Lambda (\tilde{\lambda}) \tilde{U}$, where $\Lambda (\tilde{\lambda})$ is a diagonal matrix with diagonal entries as transformed function of the sample eigenvalues and $\tilde{U}$ is the matrix of the eigenvectors. In another interesting paper \citet{Ledoit:17} proposed an estimator that shrinks the sample covariance matrix towards the identity matrix. In another paper, \citet{Karoui:15} proposed a non-parametric estimation of spectrum of eigenvalues and show that his estimator is consistent in the sense of weak convergence of distributions.

The covariance matrix estimates based on eigen-spectrum shrinkage are well-conditioned in the sense that their eigenvalues are well bounded away from zero. These estimates are based on the shrinkage of the eigenvalues and therefore invariant under some orthogonal group i.e. the shrinkage estimators shrink the eigenvalues but eigenvectors remain unchanged. In other words, the basis (eigenvector) in which the data are given is not taken advantage of and therefore the methods rely on premise that one will be able to find a good estimate in any basis. In particular, it is reasonable to believe that the basis generating the data is somewhat nice. Often this translates into the assumption that the covariance matrix has particular structure that one should be able to take advantage of. In these situations, it becomes natural to perform certain form of regularization directly on the entries of the sample covariance matrix.

Much of the recent literature focuses on two broad clases of regularized covariance matrix estimation. i) The one class relies on natural ordering among variables, where one often assumes that the variables far apart are weekly correlated and ii) the other class where there is no assumption on the natural ordering among variables. The first class includes the estimators based on banding and tapering [\citet{Bickel:3}, \citet{Cai1:6}]. These estimators are appropriate for a number of applications for ordered data (time series, spectroscopy, climate data). However for many applications including gene expression data, prior knowledge of any canonical ordering is not available and searching for all permutation of possible ordering would not be feasible. In these situations, an $\ell_1$ penalized estimator becomes more appropriate which yields a permutation-invariant estimate.

To obtain a suitable estimate which is both well-conditioned and sparse, we introduce two regularization terms: \textbf{i)} $\ell_1$ penalty for each of the off-diagonal elements of matrix and, \textbf{ii)} penalty propotional to the variance of the eigenvalues. The $\ell_1$ minimization problems are well studied in the covariance and inverse covariance matrix estimation literature [\citet{Freidman:11}, \citet{Banerjee:38}, \citet{Ravi1:24}, \citet{Bein:13}, \citet{Maurya:19} etc.]. \citet{Roth1:26} proposes an $\ell_1$ penalized log-likelihood estimator and shows that estimator is consistent in Frobenius norm at the rate of $O_P\Big(\sqrt{\{(p+s)~log~p\}/{n}}\Big)$, as both $p$ and $n$ approach to infinity. Here $s$ is the number of non-zero off-diagonal elements in the true covariance matrix. In another interesting paper \citet{Bein:13} propose an estimator of covariance matrix as penalized maximum likelihood estimator with a weighted lasso type penalty. In these optimization problems, the $\ell_1$ penalty results in sparse and a permutation-invariant estimator as compared to other $l_q, q \neq 1$ penalties. Another advantage is that the $\ell_1$ norm is a convex function which makes it suitable for large scale optimization problems. A number of fast algorithms exist in the literature for covariance and inverse covariance matrix estimation [(\citet{Freidman:11}, \citet{Roth:25}]. The eigenvalues variance penalty overcomes the over-dispersion in the sample covariance matrix  so that the estimator remains well-conditioned.

\citet{Ledoit:17} proposed an estimator of covariance matrix as a linear combination of sample covariance and identity matrix. Their estimator of covariance matrix is well-conditioned but it is not sparse. \citet{Roth:25} proposed estimator of covariance matrix based on quadratic loss function and $\ell_1$  penalty with a log-barrier on the determinant of covariance matrix. The log-determinant barrier is a valid technique to achieve positive definiteness but it is still unclear whether the iterative procedure proposed in \citet{Roth:25} actually finds the right solution to the corresponding optimization problem. In another interesting paper, \citet{Xue:31} proposed an estimator of covariance matrix as a minimizer of penalized quadratic loss function over set of positive definite matrices. In their paper, the authors solve a positive definite constrained optimization problem and establish the consistency of estimator. The resulting estimator is sparse and positive definite but whether it overcomes the over-dispersion of the eigen-spectrum of sample covariance matrix, is hard to justify. \citet{Maurya:19} proposed a joint convex penalty as function of $\ell_1$ and trace norm (defined as sum of singular values of a matrix) for inverse covariance matrix estimation based on penalized likelihood approach.

In this paper, we propose the JPEN (Joint PENalty) estimators for covariance and inverse covariance matrices estimation and derive an explicit rate of convergence in both the operator and Frobenius norm. The JPEN estimators achieves minimax rate of convergence under operator norm for the underlying class of sparse covariance and inverse covariance matrices and hence is optimal. For more details see section $\S3$. One of the major advantage of the proposed estimators is that the proposed algorithm is very fast, efficient and easily scalable to a large scale data analysis problem.

The rest of the paper is organized as following. The next section highlights some background and problem set-up for covariance and inverse covariance matrix estimation. In section 3, we describe the proposed estimators and establish their theoretical consistency. In section 4, we give an algorithm and compare its computational time with some other existing algorithms. Section 5 highlights the performance of the proposed estimators on simulated data while an application of proposed estimator to real life data is given in section 6.

\textbf{Notation:} For a matrix $M$, let $\|M\|_1 $ denote its $\ell_1$ norm defined as the sum of absolute values of the entries of $M$, $\|M\|_F$ denote its Frobenius norm, defined as the sum of square of elements of $M$, $\|M\|$ denote its operator norm (also called spectral norm), defined as the largest absolute eigenvalue of $M$, $M^{-}$ denotes matrix $M$ where all diagonal elements are set to zero, $M^{+}$ denote matrix $M$ where all off-diagonal elements are set to zero, $\sigma_i(M)$ denote the $i^{th}$ largest eigenvalue of $M$, $tr(M)$ denotes its trace, $det(M)$ denote its determinant, $\sigma_{min}(M) $ and $\sigma_{max}(M)$ denote the minimum and maximum eigenvalues of $M$, $|M|$ be its cardinality, and let $\text{sign}(M)$ be matrix of \text{sign}s of elements of $M$. For any real $x$, let $\text{sign}(x) $ denotes \text{sign} of $x$, and let $|x|$ denotes its absolute value. 

\section{Background and Problem Set-up}
Let $X=(X_1, X_2, \cdots, X_p) $ be a zero-mean p-dimensional random vector. The focus of this paper is the estimation of the covariance matrix $\Sigma:=\mathbb{E}(XX^T)$ and its inverse $\Sigma^{-1}$ from a sample of independently and identically distributed data $\{ X^{(k)} \}^{n}_{k=1}$.  In this section we provide some background and problem setup more precisely.

The choice of loss function is very crucial in any optimization problem. An optimal estimator for a particular loss function may not be optimal for another choice of loss function. Recent literature in covariance matrix and inverse covariance matrix estimation mostly focuses on estimation based on likelihood function or quadratic loss function [\citet{Freidman:11}, \citet{Banerjee:38}, \citet{Bickel:3}, \citet{Ravi1:24}, \citet{Roth:25}, \citet{Maurya:19}]. The maximum likelihood estimation requires a tractable probability distribution of observations whereas quadratic loss function does not have any such requirement and therefore fully non-parametric. The quadratic loss function is convex and due to this analytical tractability, it is a widely applicable choice for many data analysis problems. 
\subsection{Proposed Estimators}
Let $S$ be the sample covariance matrix. Consider the following optimization problem.
\begin{equation}
\hat{\Sigma}_{\lambda,\gamma}=\argmin_{\Sigma=\Sigma^T,tr(\Sigma)=tr(S)}~~\Big[  ||\Sigma-S||^2_2 +  \lambda \|{\Sigma^-}\|_1  + \gamma\sum_{i=1}^{p} \big \{\sigma_i(\Sigma)-\bar{\sigma}_{\Sigma} \big \}^2\Big],
\end{equation}
where $\bar{\sigma}_\Sigma$ is the mean of eigenvalues of $\Sigma$, $\lambda$ and $\gamma$ are some positive constants. Note that by penalty function $\|{\Sigma^-}\|_1$, we only penalize off-diagonal elements of $\Sigma$. The eigenvalues variance penalty term for eigen-spectrum shrinkage is chosen from the following points of interest: i) It is easy to interpret and ii) this choice of penalty function yields a very fast optimization algorithm. By constraint $tr(\Sigma)=tr(S)$, the total variation in $\hat{\Sigma}_{\lambda,\gamma}$ is same as that in sample covariance matrix $S$, however the eigenvalues of $\hat{\Sigma}_{\lambda,\gamma} $ are well-conditioned than those of $S$. From here onwards we suppress the dependence of $\lambda, \gamma $ on $\hat{\Sigma	}$ and denote $\hat{\Sigma}_{\lambda,\gamma} $ by $\hat{\Sigma}$. \\ \\
For $\gamma=0$, the solution to (2.1) is the standard soft-thresholding estimator for quadratic loss function and its solution is given by (see $\S4$ for derivation of this estimator):
\begin{align}
\begin{split}
\hat{\Sigma}_{ii}& =s_{ii} \\
\hat{\Sigma}_{ij}& =\text{sign}(s_{ij})\max\Big (|s_{ij}|-\frac{\lambda}{2},0\Big), ~~~~~~~~~~~~~i \neq j.
\end{split}
\end{align}
It is clear from this expression that a sufficiently large value of $\lambda$ will result in sparse  covariance matrix estimate. But estimator $\hat{\Sigma}$ of (2.2) is not necesarily positive definite [for more details here see \citet{Xue:31}]. Moreover it is hard to say whether it overcomes the over-dispersion in the sample eigenvalues. The following eigenvalue plot (Figure (2.1)) illustrates this phenomenon for a neighbourhood type (see $\S5$ for details on description of neighborhood type of covariance matrix) covariance matrix. Here we simulated random vectors from multivariate normal distribution with sample size $n=50$ and number of covariates $~p=20$.
\begin{figure}[t]
 \caption{\textit{Comparison of Eigenvalues of Covariance Matrices} }
\begin{center}
 \includegraphics[width=.8\textwidth]{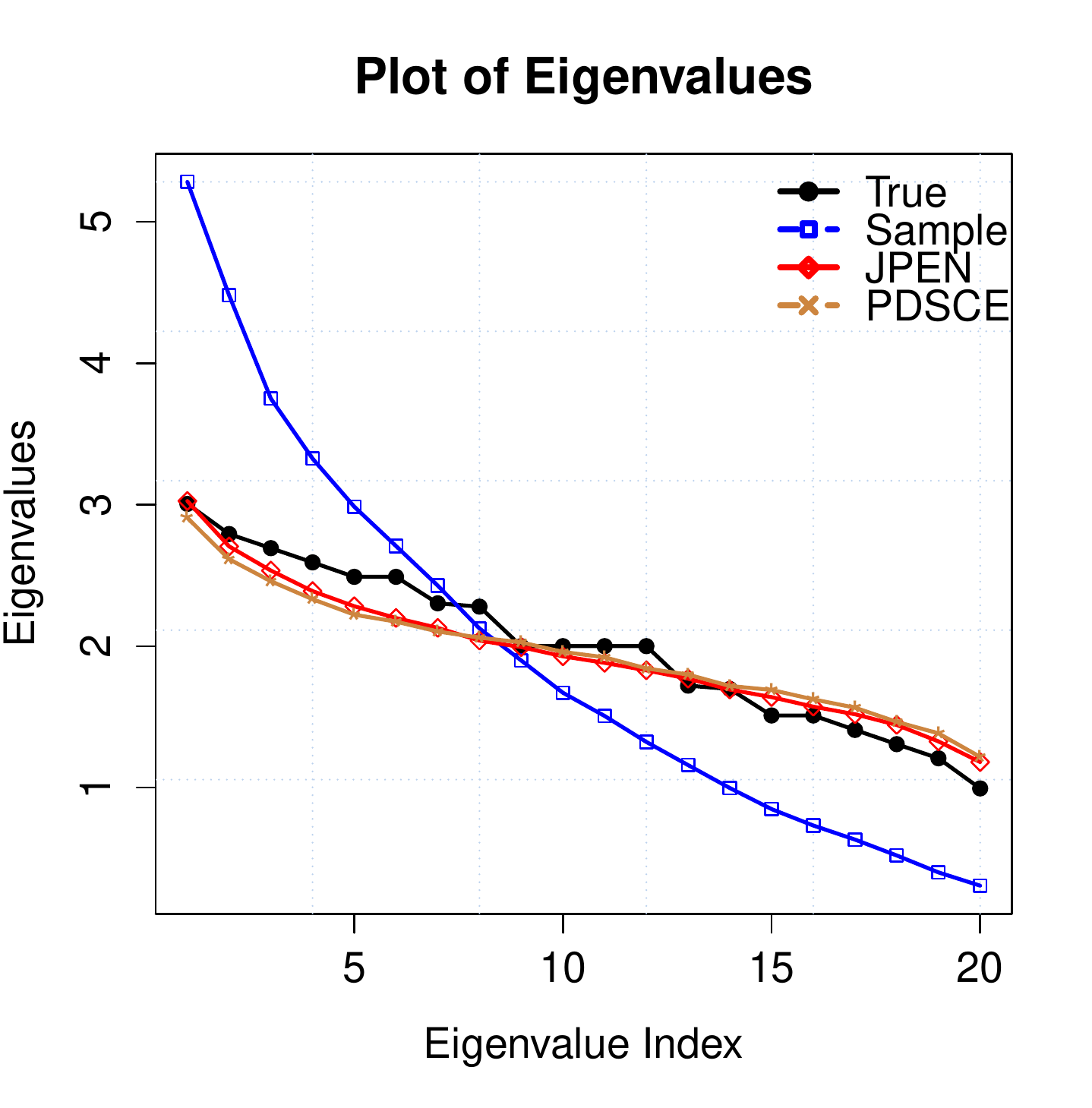}
\end{center}
\end{figure}
As is evident from Figure 2.1, eigenvalues of sample covariance matrix are over-dispersed as most of them are either too large or close to zero. Eigenvalues of the proposed Joint Penalty (JPEN) estimator and PDSCE (Positive Definite Sparse Covariance matrix Estimator (\citet{Roth1:26}) of the covariance matrix are well aligned with those of true covariance matrix. See $\S 5$ for detailed discussion. Another drawback of the estimator (2.2) is that the estimate can be negative definite.\\

As argued earlier, to overcome the over-dispersion in eigen-spectrum of sample covariance matrix, we include eigenvalues variance penalty. To illustrate its advantage, consider $\lambda=0$. After some algebra, let $\hat{\Sigma}$ be the minimizer of (2.1), then it is given by:
\begin{equation}
\hat{\Sigma}= (S+\gamma~t~I)/(1+\gamma),
\end{equation}
where $I$ is the identity matrix, and $t=\sum_{i=1}^{p}S_{ii}/p$. After some algebra, conclude that for any $\gamma>0$:
\begin{eqnarray*}
\sigma_{min} (\hat{\Sigma}) & = & \sigma_{min} (S+ \gamma~t~I)/(1+\gamma) \\
& \geq & \frac{\gamma~t}{1+\gamma}>0
\end{eqnarray*}
This means that the eigenvalues variance penalty improves $S$ to a positive definite estimator $\hat{\Sigma}$. However the estimator (2.3) is well-conditioned but need not be sparse. Sparsity can be achieved by imposing $\ell_1$ penalty on the entries of covariance matrix. Simulations have shown that, in general the minimizer of (2.1) is not positive definite for all values of $\lambda >0$ and $\gamma >0$. Here onwards we focus on correlation matrix estimation, and later generalize the method for covariance matrix estimation.\\
To achieve both well-conditioned and sparse positive definite estimator we optimize the following objective function in $R$ over specific region of values of $(\lambda, \gamma)$ which depends upon sample correlation matrix $K$ and $\lambda,\gamma$. Here the condition $tr(\Sigma)=tr(S)$ reduces to $tr(R)=p$, and $t=1$. Consider the following optimization problem:
\begin{equation}
\hat{R}_K=\argmin_{R=R^T,tr(R)=p|(\lambda, \gamma) \in \hat{\mathscr{S}}^{K}_1}~~\Big[  ||R-K||^2_F +  \lambda \|R^-\|_1  + \gamma\sum_{i=1}^{p} \big \{\sigma_i(R)-\bar{\sigma}_{R} \big \}^2\Big],
\end{equation}
where
\begin{center}
\begin{eqnarray*}
\hat{\mathscr{S}}^{K}_1& = \Big \{(\lambda,\gamma): \lambda, \gamma >0, \lambda \asymp \gamma \asymp \sqrt{\frac{log~ p}{n}}, \forall \epsilon >0,\sigma_{min}\{(K+\gamma I)-\frac{\lambda}{2}*sign(K+\gamma I)\}>\epsilon \Big \},
\end{eqnarray*}
\end{center}
and $\bar{\sigma}_{R}$ is mean of the eigenvalues of $R$. For instance when $K$ is diagonal matrix, the set $\hat{\mathscr{S}}^{K}_1$ is given by:
\begin{center}
$\hat{\mathscr{S}}^{K}_1 = \Big \{(\lambda,\gamma): \lambda, \gamma >0, \lambda \asymp \gamma \asymp \sqrt{\frac{log~ p}{n}}, \forall \epsilon >0,\lambda <2(\gamma-\epsilon) \Big \}$. 
\end{center}
The minimization in (2.4) over $R$ is for fixed $ (\lambda, \gamma) \in \hat{\mathscr{S}}^{K}_1$. The proposed estimator of covariance matrix (based on regularized correlation matrix estimator $\hat{R}_K$) is given by $\hat{\Sigma}_K=({S^+})^{1/2}\hat{R}_K({S^+})^{1/2}$, where $S^+$ is the diagonal matrix of the diagonal elements of $S$. Furthermore Lemmas 3.1 and 3.2, respectively show  that the objective function (2.4) is convex and estimator given in (2.4) is positive definite. 
\subsection{Our Contribution}
The main contributions are the following: \\
\textbf{i)} The proposed estimators are both sparse and well-conditioned simultaneously. This approach allows to take advantage of a prior structure if known on the eigenvalues of the true covariance and the inverse covariance matrices.\\
\textbf{ii)} We establish theoretical consistency of proposed estimators in both operator and Frobenius norm. The proposed JPEN estimators achieves the minimax rate of convergence in operator norm for the underlying class of sparse and well-conditioned covariance and inverse covariance matrices and therefore is optimal.\\
\textbf{iii)} The proposed algorithm is very fast, efficient and easily scalable to large scale optimization problems.
\section{Analysis of JPEN Method}
\textbf{Def:} A random vector $X$ is said to have sub-gaussian distribution if for each $t \ge 0$ and $y \in \mathbb{R}^p $ with $\|y\|_2=1$, there exist $0< \tau < \infty $ such that
\begin{equation}
\mathbb{P}\{|y^T(X-\mathbb{E}(X))|>t\} \le e^{-t^2/2\tau}
\end{equation}
Although the JPEN estimators exists for any finite $2 \le n<p<\infty$, for theoretical consistency in operator norm we require $s~log~p=o(n)$ and for Frobenus norm we require $(p+s) ~log~p=o(n)$ where $s$ is the upper bound on the number of non-zero off-diagonal entries in true covariance matrix. For more details, see the remark after Theorem 3.1. \\
\subsection{Covariance Matrix Estimation}
We make the following assumptions about the true covariance matrix $\Sigma_0$.\\
\textbf{A0.} Let $X:=(X_1,X_2,\cdots,X_p)$ be a mean zero vector with covariance matrix $\Sigma_0$ such that each $X_i/ \sqrt{\Sigma_{0ii}}$ has subgaussian distribution with parameter $\tau$ as defined in (3.1).\\
\textbf{A1.} With $ E=\{(i,j): \Sigma_{0ij} \neq 0, i \neq j \}, $ the $|E| \le s $ for some positive integer $s$. \\
\textbf{A2.} There exists a finite positive real number $\bar{k} >0$ such that $ 1/\bar{k} \le \sigma_{min}(\Sigma_0) \le \sigma_{max}(\Sigma_0) \le \bar{k}$.\\

Assumption A2 guarantees that the true covariance matrix $\Sigma_0$ is well-conditioned (i.e. all the eigenvalues are finite and positive). A well-conditioned means that [\citet{Ledoit:17})] inverting the matrix does not explode the estimation error.  Assumption A1 is more of a definition which says that the number of non-zero off diagonal elements are bounded by some positive integer. Theorem 3.1 gives the rate of convergence of the proposed correlation based covariance matrix estimator (2.4). The following Lemmas show that optimization problem in (2.4) is convex and the proposed JPEN estimator (2.4) is positive definite.
\begin{Lemma}
The optimization problem in (2.4) is convex.
\end{Lemma}
\begin{Lemma}
The estimator given by (2.4) is positive definite for any $2 \le n < \infty $ and $p <\infty$.
\end{Lemma}
\begin{Theorem}
Let $(\lambda, \gamma) \in \hat{\mathscr{S}}^{K}_1$ and $\hat{\Sigma}_K$ be as defined in (2.4). Under Assumptions A0, A1, A2,
\begin{equation}
\|\hat{R}_K-R_0\|_F=O_P \Big ( \sqrt{\frac{s~ log~p}{n}} \Big ) ~~~ \text{and} ~~~\|\hat{\Sigma}_K-\Sigma_0\|=O_P \Big ( \sqrt{\frac{(s+1) log~p}{n}} \Big ),
\end{equation}
where $R_0$ is true correlation matrix.
\end{Theorem}
\textbf{Remark: 1.} The JPEN estimator $\hat{\Sigma}_K$ is minimiax optimal under the operator norm. In (\citet{Cai2:40}), the authors obtain the minimax rate of convergence in the operator norm of their covariance matrix estimator for the particular construction of parameter space $\mathscr{H}_0(c_{n,p}):=\Big \{ \Sigma : max_{1 \le i \le p}\sum_{i=1}^{p}I\{\sigma_{ij}\neq 0\} \leq c_{n,p} \Big \}$. They show that this rate in operator norm is $c_{n,p} \sqrt{log~p/n}$ which is same as that of $\hat{\Sigma}_K$ for $1 \leq c_{n,p}=\sqrt{s}$. \\

{\bf 2.} \citet{Bickel1:4} proved that under the assumption of $\sum_{j=1}|\sigma_{ij}|^q \leq c_0(p)$ for some $ 0 \leq q \leq 1$, the hard thresholding estimator of the sample covariance matrix for tuning parameter $\lambda \asymp \sqrt{(log~p)/n}$ is consistent in operator norm at a rate no worse than $ O_P\Big ( c_0(p) \sqrt{p}(\frac{log ~p}{n})^{(1-q)/2} \Big ) $ where $c_0(p)$ is the upper bound on the number of non-zero elements in each row. Here the truly sparse case corresponds to $q=0$. The rate of convergence of $\hat{\Sigma}_K$ is same as that of \citet{Bickel1:4} except in the following cases:\\
{\bf Case (i)} The covariance matrix has all off diagonal elements zero except last row which has $ \sqrt{p}$ non-zero elements. Then $c_0(p)=\sqrt{p}$ and $ \sqrt{s}=\sqrt{2~\sqrt{p}-1}$. The opeartor norm rate of convergence for JPEN estimator is $O_P \Big ( \sqrt{\sqrt{p}~(log~p)/n} \Big )$ where as rate of Bickel and Levina's estimator is $O_P \Big (\sqrt{p~(log~p)/n} \Big )$. \\
{\bf Case (ii)} When the true covariance matrix is tridiagonal, we have $c_0(p)=2$ and $s=2p-2$, the JPEN estimator has rate of $\sqrt{p~log~p/n}$ whereas the Bickel and Levina's estimator has rate of $\sqrt{log~p/n}$. \\
For the case $\sqrt{s} \asymp c_0(p)$ and JPEN has the same rate of convergence as that of Bickel and Levina's estimator.\\

\textbf{3.} The operator norm rate of convergence is much faster than Frobenius norm. This is due to the fact that Frobenius norm convergence is in terms of all eigenvalues of the covariance matrix whereas the operator norm gives the convergence of the estimators in terms of the largest eigenvalue.\\

\textbf{4.} Our proposed estimator is applicable to estimate any non-negative definite covariance matrix.\\

Note that the estimator $\hat{\Sigma}_K$ is obtained by regularization of sample correlation matrix in (2.4). In some application it is desirable to directly regularize the sample covariance matrix. The JPEN estimator of the covariance matrix based on regularization of sample covariance matrix is obtained by solving the following optimization problem:
\begin{equation}
\hat{\Sigma}_S=\argmin_{\Sigma=\Sigma^T,tr(\Sigma)=tr(S)|(\lambda, \gamma) \in \hat{\mathscr{S}}^{S}_1}~~\Big[  ||\Sigma-S||^2_F +  \lambda \|\Sigma^-\|_1  + \gamma\sum_{i=1}^{p} \{\sigma_i(\Sigma)-\bar{\sigma}_{\Sigma}\}^2\Big],
\end{equation}
where
\begin{eqnarray*}
\hat{\mathscr{S}}^S_1& = \Big \{(\lambda,\gamma): \lambda,\gamma >0, \lambda \asymp \gamma \asymp \sqrt{\frac{log~ p}{n}},  \forall \epsilon >0, \sigma_{min}\{(S+\gamma t I)-\frac{\lambda}{2}*sign(S+\gamma t I)\}>\epsilon \},
\end{eqnarray*}
and $S$ is sample covariance matrix. The minimization in (3.3) over $\Sigma$ is for fixed $ (\lambda, \gamma) \in \hat{\mathscr{S}}^{S}_1$. The estimator $\hat{\Sigma}_S$ is positive definite and well-conditioned. Theorem 3.2 gives the rate of convergence of the estimator $\hat{\Sigma}_S$ in Frobenius norm.
\begin{Theorem}
Let $(\lambda, \gamma) \in \hat{\mathscr{S}}^S_1$, and let $\hat{\Sigma}_S$ be as defined in (3.3). Under Assumptions A0, A1, A2,
\begin{equation}
\|\hat{\Sigma}_S-\Sigma_0\|_F=O_P \Big ( \sqrt{\frac{(s+p) log~p}{n}} \Big )
\end{equation}
\end{Theorem}
As noted in \citet{Roth1:26} the worst part of convergence here comes from estimating the diagonal entries. 
\subsubsection{Weighted JPEN Estimator for the Covariance Matrix Estimation}
A modification of estimator $\hat{R}_{K}$ is obtained by adding positive weights to the term $(\sigma_i(R)-\bar{\sigma}_R)^2$. This leads to weighted eigenvalues variance penalty with larger weights amounting to greater shrinkage towards the center and vice versa. Note that the choice of the weights allows one to use any prior structure of the eigenvalues (if known) in estimating the covariance matrix. The weighted JPEN correlation matrix estimator $\hat{R}_A$ is given by :
\begin{equation}
\hat{R}_A=\argmin_{R=R^T,tr(R)=p|(\lambda, \gamma) \in \hat{\mathscr{S}}^{K,A}_1}~~\Big[  ||R-K||^2_F +  \lambda \|R^-\|_1  + \gamma\sum_{i=1}^{p} a_i\{\sigma_i(R)-\bar{\sigma}_R\}^2\Big],
\end{equation}
where
\begin{eqnarray*}
\hat{\mathscr{S}}^{K,A}_1& = \Big \{(\lambda,\gamma): \lambda \asymp \gamma \asymp \sqrt{\frac{log~p}{n}}, \lambda \le \frac{(2~\sigma_{min}(K))(1+\gamma ~max(A_{ii})^{-1})}{(1+\gamma~ min(A_{ii}))^{-1}p}+ \frac{\gamma ~min(A_{ii})}{p} \Big \},
\end{eqnarray*}
and $A=\text{diag}(A_{11},A_{22},\cdots A_{pp})$ with $A_{ii}=a_i$. The proposed covariance matrix estimator is $\hat{\Sigma}_{K,A}=(S^{+})^{1/2}\hat{R}_A (S^{+})^{1/2}$. The optimization problem in (3.5) is convex and yields a positive definite estimator for each $(\lambda,\gamma) \in \hat{\mathscr{S}}^{K,A}_1$. A simple excercise shows that the estimator $\hat{\Sigma}_{K,A}$ has same rate of convergence as that of $\hat{\Sigma}_{S}$.

\subsection{Estimation of Inverse Covariance Matrix}
We extend the JPEN approach to estimate a well-conditioned and sparse inverse covariance matrix. Similar to the covariance matrix estimation, we first propose an estimator for inverse covariance matrix based on regularized inverse correlation matrix and discuss its rate of convergence in Frobenious and operator norm. \\

\textbf{Notation:} We shall use $Z$ and $\Omega$ for inverse correlation and inverse covariance matrix respectively. \\
\textbf{Assumptions:} We make the following assumptions about the true inverse covariance matrix $\Omega_0$. Let $\Sigma_0=\Omega_0^{-1}$.\\
\textbf{B0.} Same as the assumption $A0$.\\
\textbf{B1.} With $ H=\{(i,j): \Omega_{0ij} \neq 0, i \neq j \}$, the $|H| \le s $, for some positive integer $s$. \\
\textbf{B2.} There exist $ 0< \bar{k} < \infty $ large enough such that $ (1/{\bar{k}}) \le \sigma_{min}(\Omega_0) \le \sigma_{max}(\Omega_0) \le \bar{k}$.\\

Let $\hat{R}_K$ be a JPEN estimator for the true correlation matrix. By Lemma 3.2, $\hat{R}_K$ is positive definite. Define the JPEN estimator of inverse correlation matrix as the solution to the following optimization problem,
\begin{equation}
\hat{Z}_{K}=\argmin_{Z=Z^T,tr(Z)=tr(\hat{R}_K^{-1})|(\lambda,\gamma) \in \hat{\mathscr{S}}^{{K}}_2}\Big [ \|Z-\hat{R}_K^{-1} \|^2~+~\lambda\|Z^-\|_1~+~\gamma\sum_{i=1}^{p} \{\sigma_i(Z)- \bar{\sigma}(Z)\}^2 \Big ]
\end{equation}
where
\begin{center}
\begin{align*}
\hat{\mathscr{S}}^{{K}}_2& = \Big \{(\lambda,\gamma): \lambda,\gamma >0, \lambda \asymp \gamma \asymp \sqrt{\frac{log~ p}{n}}, \forall \epsilon >0, \\
& ~~~~~~~~~~~~ \sigma_{min}\{(\hat{R}_K^{-1}+\gamma t_1 I)-\frac{\lambda}{2}*sign(\hat{R}_K^{-1}+\gamma t_1 I)\}>\epsilon \Big \},
\end{align*}
\end{center}
and $t_1$ is average of the diagonal elements of $\hat{R}_K^{-1}$. The minimization in (3.6) over $Z$ is for fixed $ (\lambda, \gamma) \in \hat{\mathscr{S}}^{{K}}_2$. The proposed JPEN  estimator of inverse covariance matrix (based on regularized inverse correlation matrix estimator $\hat{Z}_{K}$) is given by $\hat{\Omega}_{{K}}=(S^+)^{-1/2}{\hat{Z}}_{{K}}(S^+)^{-1/2}$, where $S^+$ is a diagonal matrix of the diagonal elements of $S$. Moreover (3.6) is a convex optimization problem and $\hat{Z}_K$ is positive definite. \\

Next we state the consistency of estimators $\hat{Z}_{{K}}$ and $\hat{\Omega}_{{K}}$.
\begin{Theorem}
Under Assumptions B0, B1, B2 and for $(\lambda,\gamma) \in \hat{\mathscr{S}}^{{K}}_{2}$,
\begin{equation}
\|\hat{Z}_{{K}}-R_0^{-1}\|_F=O_P \Big ( \sqrt{\frac{s~log~p}{n}} \Big )  ~~~\text{and} ~~~\|\hat{\Omega}_{{K}}-\Omega_0\|=O_P \Big ( \sqrt{\frac{(s+1)~log~p}{n}} \Big )
\end{equation}
where $R_0^{-1}$ is the inverse of true correlation matrix.
\end{Theorem}
\textbf{Remark:1.} Note that the JPEN estimator $\hat{\Omega}_{{K}}$ achieves minimax rate of convergence for the class of covariance matrices satisfying assumption $B0$, $B1$, and $B2$ and therefore optimal. The similar rate is obtained in \citet{Cai2:40} for their class of sparse inverse covariance matrices.\\

Next we give another estimate of inverse covariance matrix based on $\hat{\Sigma}_{S}$. Consider the following optimization problem:
\begin{equation}
\hat{\Omega}_S=\argmin_{\Omega=\Omega^T,tr(\Omega)=tr(\hat{\Sigma}_{S}^{-1})|(\lambda, \gamma) \in \hat{\mathscr{S}}^{{S}}_2}~~\Big[  ||\Omega- \hat{\Sigma}_{S}^{-1}||^2_F +  \lambda \|\Omega^-\|_1  + \gamma\sum_{i=1}^{p} \{\sigma_i(\Omega)-\bar{\sigma}_{\Omega}\}^2\Big],
\end{equation}
where
\begin{center}
\begin{align*}
\hat{\mathscr{S}}^{{S}}_2&= \Big \{(\lambda,\gamma): \lambda,\gamma >0, \lambda \asymp \gamma \asymp \sqrt{\frac{log~ p}{n}}, ~\forall \epsilon >0, \\
& ~~~~~~~~~~~~\sigma_{min}\{(\hat{\Sigma}_S^{-1}+\gamma ~t_2~I)-\frac{\lambda}{2}*sign(\hat{\Sigma}_S^{-1}+\gamma t_2 I)\}>\epsilon \Big \},
\end{align*}
\end{center}
and $t_2$ is average of the diagonal elements of $\hat{\Sigma}_S $. The minimization in (3.8) over $\Omega$ is for fixed $ (\lambda, \gamma) \in \hat{\mathscr{S}}^{S}_2$. The estimator in (3.8) is positive definite and well-conditioned. The consistency result of the estimator $\hat{\Omega}_S$ is given in following theorem.
\begin{Theorem}
Let $(\lambda, \gamma) \in \hat{\mathscr{S}}^{S}_2$ and let $\hat{\Omega}_S$ be as defined in (3.8). Under Assumptions B0, B1, and B2,
\begin{equation}
\|\hat{\Omega}_S-\Omega_0\|_F=O_P \Big ( \sqrt{\frac{(s+p) log~p}{n}} \Big ).
\end{equation}
\end{Theorem}
\subsubsection{Weighted JPEN Estimator for The Inverse Covariance Matrix}
Similar to weighted JPEN covariance matrix estimator $\hat{\Sigma}_{K,A}$, a weighted JPEN estimator of the inverse covariance matrix is obtained by adding positive weights $a_i $ to the term $(\sigma_i(Z)-1)^2$ in (3.8). The weighted JPEN estimator is $\hat{\Omega}_{{K},A}:=({S^{+}})^{-1/2}\hat{Z}_A ({S^{+}})^{-1/2}$, where
\begin{equation}
\hat{Z}_A=\argmin_{Z=Z^T,tr(Z)=tr(\hat{R}_K^{-1})|(\lambda, \gamma) \in \hat{\mathscr{S}}^{K,A}_2}~~\Big[  ||Z-\hat{R}_K^{-1}||^2_F +  \lambda \|Z^-\|_1  + \gamma\sum_{i=1}^{p} a_i\{\sigma_i(Z)-1\}^2\Big],
\end{equation}
with
\begin{eqnarray*}
\hat{\mathscr{S}}^{{K},A}_2& = \Big \{(\lambda,\gamma): \lambda \asymp \gamma \asymp \sqrt{\frac{log~p}{n}}, \lambda \le \frac{(2~\sigma_{min}({R}_K^{-1}))(1+\gamma t_1 max(A_{ii})^{-1})}{(1+\gamma~ min(A_{ii})^{-1}p}+ \frac{\gamma min(A_{ii})}{p} \Big \},
\end{eqnarray*}
and $A=\text{diag}(A_{11},A_{22},\cdots A_{pp})$ with $A_{ii}=a_i$. The optimization problem in (3.10) is convex and yields a positive definite estimator for $(\lambda,\gamma) \in \hat{\mathscr{S}}^{K,A}_2$. A simple excercise shows that the estimator $\hat{Z}_{A}$ has similar rate of convergence as that of $\hat{Z}_K$.

\section{An Algorithm} 
\subsection{Covariance Matrix Estimation:}
The optimization problem (2.4) can be written as:\\
\begin{eqnarray}
\hat{R}_K=\argmin_{R=R^T|(\lambda, \gamma) \in \hat{\mathscr{S}}^{K}_1 }~f(R),
\end{eqnarray}
where
\begin{eqnarray*}
f(R)=||R-K||^2_F +  \lambda \|{R}^-\|_1  + \gamma\sum_{i=1}^{p} \{\sigma_i(R)-\bar{\sigma}(R)\}^2.
\end{eqnarray*}
Note that $\sum_{i=1}^{p} \{\sigma_i(R)-\bar{\sigma}(R)\}^2=tr(R^2)-2~tr(R)+p$, where we have used the constraint $tr(R)=p$. Therefore,
\begin{eqnarray*}
f(R) & = & \|R-K\|_F^2+\lambda \|{R}^-\|_1 +\gamma~ tr(R^2)-2~ \gamma~tr(R) +p\\
&  = & tr(R^2 (1+\gamma))-2tr\{R(K+\gamma I)\}+ tr(K^TK) +\lambda ~\|{R}^-\|_1 +p\\
&  = & (1+\gamma)\{tr(R^2)-2/(1+\gamma) tr\{R(K+\gamma I)\}+ (1/(1+\gamma))tr(K^TK)\} \\
&   &~~~~+\lambda ~\|{R}^-\|_1 +p \\
&  = & (1+\gamma)\{\|R- (K+\gamma I)/(1+\gamma)\|^2_F+ (1/(1+\gamma))tr(K^TK)\} \\
 &   & ~~~~+\lambda ~\|{R}^-\|_1 +p.
\end{eqnarray*}
The solution of (4.1) is soft thresholding estimator and it is given by:
\begin{eqnarray}
\hat{R}_K =\frac{1}{1+\gamma}~\text{sign}(K)*\text{pmax}\{\text{abs}(K+\gamma ~I)-\frac{\lambda}{2},0\}
\end{eqnarray}
with $(\hat{R}_{K})_{ii}=(K_{ii}+\gamma)/(1+\gamma)$, $pmax(A,b)_{ij}:=max(A_{ij},b)$ is elementwise max function for each entry of the matrix $A$. Note that for each $(\lambda,\gamma) \in \hat{\mathscr{S}}^{K}_1$, $\hat{R}_{K}$ is positive definite. \\ \\
{\bf Choice of $\lambda$ and $\gamma$:}
For a given value of $\gamma$, we can find the value of $\lambda $ satisfying:
\begin{eqnarray}
\sigma_{min}\{(K+\gamma I)-\frac{\lambda}{2}*sign(K+\gamma I)\}>0
\end{eqnarray}
which can be simplified to 
\begin{align*}
\lambda < \frac{\sigma_{min}(K+\gamma I)}{C_{12}~\sigma_{max}(\text{sign}(K))}. 
\end{align*}
For some $C_{12} \ge 0.5$. Such choice of $(\lambda,\gamma)\in \hat{\mathscr{S}}_1^{{K}}$, and the estimator $\hat{R}_{K}$ is positive definite. Smaller values of $C_{12}$ yeild a solution which is more sparse but may not be positive definite. \\ \\
{\bf Choice of weight matrix A:}
For optimization problem in (3.5), the weights are chosen in following way:\\
Let $\mathscr{E}$ be the set of sorted diagonal elements of the sample covariance matrix $S$.\\
\textbf{i)} Let $k$ be largest index of $\mathscr{E}$ such that $k^{th}$ elements of $\mathscr{E}$ is less than $1$. For $i \le k, ~a_i=\mathscr{E}_{i}$. For $k<i \le p,~ a_i=1/\mathscr{E}_{i}.$\\ 
\textbf{ii)} $A=\text{diag}(a_1,a_2,\cdots,a_p),~\text{where}~ a_j=a_j/\sum_{i=1}^p a_i.$
Such choice of weights allows more shrinkage of extreme sample eigenvalues than the ones in center of eigen-spectrum. 
\subsection{Inverse Covariance Matrix Estimation:}
To get an expression of inverse covariance matrix estimate, we replace $K$ by $\hat{R}_K^{-1}$ in (4.2), where $\hat{R}_K$ is a JPEN estimator of correlation matrix. We chose $(\lambda,\gamma) \in \hat{\mathscr{S}}^{{K}}_2$. For a given $\gamma$, we chose $\lambda>0$ satisfying:\\
\begin{eqnarray}
\sigma_{min}\{(\hat{R}_K^{-1}+\gamma t_1 I)-\frac{\lambda}{2}*sign(\hat{R}_K^{-1}+\gamma t_1 I)\}>0
\end{eqnarray}
which can be simplified to
\begin{align*}
\lambda < \frac{\sigma_{min}(\hat{R}_K^{-1}+\gamma t_1 I)}{C_{12}~\sigma_{max}(\text{sign}(\hat{R}_K^{-1}))}. 
\end{align*}

\subsection{Computational Complexity}
The JPEN estimator $\hat{\Sigma}_{K}$ has computational complexity of $O(p^2)$ as there are at most $3p^2$ multiplications for computing the estimator $\hat{\Sigma}_{K}$. The other existing algorithm Glasso (\cite{Freidman:11}), PDSCE (\cite{Roth1:26}) have computational complexity of $O(p^3)$. We compare the computational timing of our algorithm to some other existing algorithms Glasso (\citet{Freidman:11}), PDSCE (\citet{Roth1:26}). The exact timing of these algorithm also depends upon the implementation, platform etc. (we did our computations in $R$ on a AMD 2.8GHz processor). Following the approach \cite{Bickel1:4}, the optimal tuning parameter $(\lambda,\gamma)$ was obtained by minimizing the $5-$fold cross validation error 
$$(1/5) \sum_{i=1}^5\|\hat{\Sigma}_i^v-\Sigma_i^{-v}\|_1,$$ where $\hat{\Sigma}_i^v$ is JPEN estimate of the covariance matrix \begin{wrapfigure}{H}{0.5\textwidth}
\vspace{-5pt} \begin{center}
    \includegraphics[width=0.40\textwidth]{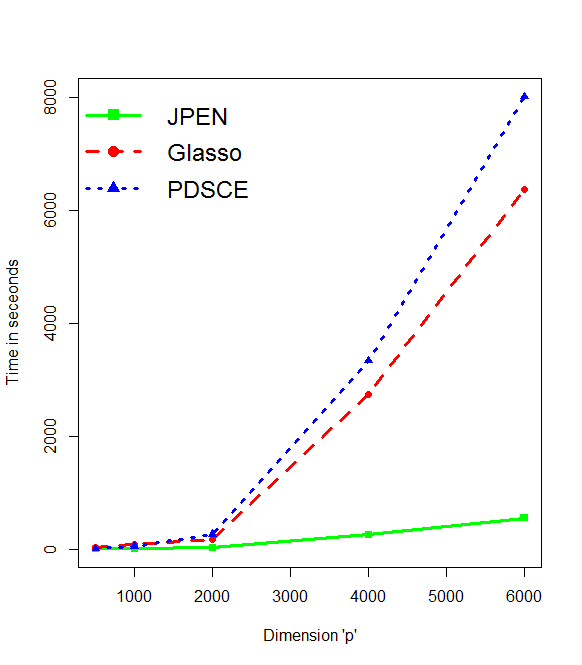}
 \end{center}
\vspace{-5pt} \caption{\textit{Timing comparison of JPEN, Glasso, and PDSCE.} }
\end{wrapfigure}
based on $v=4n/5$ observations, $\Sigma_i^{-v}$ is the sample covariance matrix using $(n/5)$ observations. 
Figure 4.1 illustrates the total computational time taken to estimate the covariance matrix by $Glasso,~PDSCE$ and $JPEN$ algorithms for different values of $p$ for Toeplitz type of covariance matrix on log-log scale (see section $\S 5$ for Toeplitz type of covariance matrix). Although the proposed method requires optimization over a grid of values of $(\lambda,\gamma) \in \hat{\mathscr{S}}^{K}_1$, our algorithm is very fast and easily scalable to large scale data analysis problems.
\section{Simulation Results}
\paragraph{•}
We compare the performance of the proposed method to other existing methods on simulated data for five types of covariance and inverse covariance matrices. \\

\textbf{(i) Hub Graph:} Here the rows/columns of $\Sigma_0$ are partitioned into J equally-sized disjoint groups: $ \{ V_1 \cup  V_2~  \cup, ..., \cup~ V_J\} = \{1,2,...,p\},$ each group is associated with a $\bf pivotal$ row k. Let size $|V_1| = s$. We set $\sigma_{0i,j}=\sigma_{0j,i}=\rho $ for $ i \in  V_k $ and $\sigma_{0i,j}=\sigma_{0j,i}=0$ otherwise. In our experiment, $J=[p/s], k =1,s+1, 2s+1,...,$ and we always take $\rho= 1/(s + 1)$ with J = 20. \\

\textbf{(ii) \bf Neighborhood Graph:} We first uniformly sample $(y_1,y_2,...,y_n)$ from a unit square. We then set $\sigma_{0i,j}=\sigma_{0j,i}=\rho $ with probability ${(\sqrt{2\pi})}^{-1}exp( -4\|y_i-y_j\|^2)$. The remaining entries of $\Sigma_0$ are set to be zero. The number of nonzero off-diagonal elements of each row or column is restricted to be smaller than $[1/\rho]$ where $\rho$ is set to be 0.245. \\

\textbf{(iii) \bf Toeplitz ~ Matrix:} We set $\sigma_{0i,j}=2~$for$ ~i=j;~\sigma_{0i,j}=|0.75|^{|i-j|}~$ for $|i-j|=1,2;$ and $~\sigma_{0i,j}=0~$ otherwise.  \\

\textbf{(iv) \bf Block ~Diagonal~ Matrix:} In this setting $\Sigma_0$ is a block diagonal matrix with varying block size. For $p=500$ number of blocks is 4 and for $p=1000$ the number of blocks is 6. Each block of covariance matrix is taken to be Toeplitz type matrix as in case (iii).\\

\textbf{(v) Cov-I~ type Matrix:} In this setting, we first simulate a random sample $(y_1,y_2,...,y_p)$ from standard normal distribution. Let $x_i=|y_i|^{3/2}*(1+1/p^{1+log(1+1/p^2)})$. Next we generate multivariate normal random vectors $\mathscr{Z}=(z_1,z_2,...,z_{5p})$ with mean vector zero and identity covariance matrix. Let $U$ be eigenvector corresponding to sample covariance matrix of $\mathscr{Z}$. We take $\Sigma_0=UDU'$, where $D=\text{diag}(x_1,x_2,....x_p)$. This is not a sparse setting but the covariance matrix has most of eigenvalues close to zero and hence allows us to compare the performance of various methods in a setting where most of eigenvalues are close to zero and widely spread as compared to structured covariance matrices in \textbf{(i)-(iv)}.
\\
\begin{table}[t]\scriptsize									
\centering									
\caption{Covariance Matrix Estimation}		
\vspace{3mm}						
\begin{tabular}{|r|p{1.6cm}|p{1.6cm}|p{1.6cm}|p{1.6cm}|}									
\hline							
\multicolumn{5}{ |c| }{Block type covariance matrix} \\									
\hline									
{} & \multicolumn{2}{c|}{n=50} & \multicolumn{2}{c|}{n=100}\\									
\hline									
& p=500     & p=1000  & p=500   & p=1000 \\									
\midrule									
Ledoit-Wolf	&	1.54(0.102)	&	2.96(0.0903)	&	4.271(0.0394)	&	2.18(0.11)	\\
Glasso	&	0.322(0.0235)	&	3.618(0.073)	&	0.227(0.098)	&	2.601(0.028)	\\
PDSCE	&	3.622(0.231)	&	4.968(0.017)	&	1.806(0.21)	&	2.15(0.01)	\\
BLThresh	&	2.747(0.093)	&	3.131(0.122)	&	0.887(0.04)	&	0.95(0.03) 	\\
JPEN	&	2.378(0.138)	&	3.203(0.144)	&	1.124(0.088)	&	2.879(0.011)	\\
\bottomrule									
									
\multicolumn{5}{ |c| }{Hub type covariance matrix} \\									
\hline									
{} & \multicolumn{2}{c|}{n=50} & \multicolumn{2}{c|}{n=100}\\									
\hline									
& p=500     & p=1000  & p=500   & p=1000 \\									
\midrule									
Ledoit-Wolf	&	2.13(0.103)	&	2.43(0.043) 	&	1.07(0.165)	&	3.47(0.0477)	\\
Glasso	&	0.511(0.047)	&	0.551(0.005)	&	0.325(0.053)	&	0.419(0.003)	\\
PDSCE	&	0.735(0.106)	&	0.686(0.006)	&	0.36(0.035)	&	0.448(0.002)	\\
BLThresh	&	1.782(0.047)	&	2.389(0.036) 	&	0.875(0.102)&		1.82(0.027)	\\
JPEN	&	0.732(0.111)	&	0.688(0.006)	&	0.356(0.058)	&	0.38(0.007)	\\
\bottomrule									
\multicolumn{5}{ |c| }{Neighborhood type covariance matrix} \\									
\hline									
{} & \multicolumn{2}{c|}{n=50} & \multicolumn{2}{c|}{n=100}\\									
\hline									
& p=500     & p=1000  & p=500   & p=1000 \\									
\midrule									
Ledoit-Wolf	&	1.36(0.054)	&	2.89(0.028)	&	1.1(0.0331)	&	2.32(0.0262)\\	
Glasso	&	0.608(0.054)	&	0.63(0.005)	&	0.428(0.047)	&	0.419(0.038)	\\
PDSCE	&	0.373(0.085)	&	0.468(0.007)	&	0.11(0.056)	&	0.175(0.005)	\\
BLThresh	&	1.526(0.074) 	&	2.902(0.033)	&	0.870(0.028)&		1.7(0.026)\\	
JPEN	&	0.454(0.0423)	&	0.501(0.018)	&	0.086(0.045)	&	0.169(0.003)	\\
\bottomrule									
\multicolumn{5}{ |c| }{Toeplitz type covariance matrix} \\									
\hline									
{} & \multicolumn{2}{c|}{n=50} & \multicolumn{2}{c|}{n=100}\\									
\hline									
& p=500     & p=1000  & p=500   & p=1000 \\									
\midrule									
Ledoit-Wolf	&	1.526(0.074)	&	  2.902(0.033)		&	1.967(0.041)	&	2.344(0.028)\\
Glasso	&	2.351(0.156)	&	3.58(0.079)		&	1.78(0.087)	&	2.626(0.019) \\
PDSCE	&	3.108(0.449)	&	5.027(0.016)		&	0.795(0.076)	&	2.019(0.01) \\
BLThresh	&	0.858(0.040)	&	1.206(0.059)		&	0.703(0.039)&		1.293(0.018)\\
JPEN	&	2.517(0.214)	&	3.205(0.16)		&	1.182(0.084)	&	2.919(0.011) \\
\bottomrule									
\multicolumn{5}{ |c| }{Cov-I type covariance matrix} \\									
\hline									
{} & \multicolumn{2}{c|}{n=50} & \multicolumn{2}{c|}{n=100}\\									
\hline									
& p=500     & p=1000  & p=500   & p=1000 \\									
\midrule									
Ledoit-Wolf	&	33.2(0.04)	&	36.7(0.03)	&	36.2(0.03)	&	48.0(0.03)	\\
Glasso	&	15.4(0.25)	&	16.1(0.4)	&	14.0(0.03)	&	14.9(0.02)	\\
PDSCE	&	16.5(0.05)	&	16.33(0.04)	&	16.9(0.03)	&	17.5(0.02)	\\
BLThresh	&	15.7(0.04)	&	17.1(0.03)	&	13.4(0.02)	&	17.5(0.02)	\\
JPEN	&	7.1(0.042)	&	11.5(0.07)	&	8.4(0.042)	&	7.8(0.034)	\\
\bottomrule									
\end{tabular}%
\label{tab:addlabel}%
\end{table}%

We chose similar structure of $\Omega_0$ for simulations. For all these choices of covariance and inverse covariance matrices, we generate random vectors from multivariate normal distribution with varying $n$ and $p$. 
We chose $n=50,100$ and $p=500,1000$. We compare the performance of proposed covariance matrix estimator $\hat{\Sigma}_K$ to graphical lasso [\citet{Freidman:11}], PDSC Estimate [\citet{Roth1:26}], Bickel and Levina's thresholding estimator (BLThresh) [\cite{Bickel1:4}] and Ledoit-Wolf [\cite{Ledoit:17}] estimate of covariance matrix. The JPEN estimate $\hat{\Sigma}_{K}$ was computed using R software(version 3.0.2). The graphical lasso estimate of the covariance matrix was computed using R package ``glasso" (http://statweb.stanford.edu/ tibs/glasso/).  
\begin{table}[t]\scriptsize									
\centering									
\caption{Inverse Covariance Matrix Estimation}		
\vspace{3mm}						
\begin{tabular}{|r|p{1.6cm}|p{1.6cm}|p{1.6cm}|p{1.6cm}|}									
\hline
\multicolumn{5}{ |c| }{Block type covariance matrix} \\									
\hline									
{} & \multicolumn{2}{c|}{n=50} & \multicolumn{2}{c|}{n=100}\\									
\hline									
& p=500     & p=1000  & p=500   & p=1000 \\									
\midrule									
Glasoo	&	4.144(0.523)	&	1.202(0.042)	&	0.168(0.136)	&	1.524(0.028)	\\
PDSCE	&	1.355(0.497)	&	1.201(0.044)	&	0.516(0.196)	&	0.558(0.032)	\\
CLIME	&	4.24(0.23)	&	6.56(0.25)	&	6.88(0.802)	&	10.64(0.822)	\\
JPEN	&	1.248(0.33)	&	1.106(0.029)	&	0.562(0.183)	&	0.607(0.03)	\\
\bottomrule									
\multicolumn{5}{ |c| }{Hub type covariance matrix} \\									
\hline									
{} & \multicolumn{2}{c|}{n=50} & \multicolumn{2}{c|}{n=100}\\									
\hline									
& p=500     & p=1000  & p=500   & p=1000 \\									
\midrule									
Glasoo	&	1.122(0.082)	&	0.805(0.007)	&	0.07(0.038)	&	0.285(0.004)	\\
PDSCE	&	0.717(0.108)	&	0.702(0.007)	&	0.358(0.046)	&	0.356(0.005)	\\
CLIME	&	10.5(0.329)	&	10.6(0.219)	&	6.98(0.237)	&	10.8(0.243)	\\
JPEN	&	0.684(0.051)	&	0.669(0.003)	&	0.34(0.024)	&	0.337(0.002)	\\
\bottomrule									
\multicolumn{5}{ |c| }{Neighborhood type covariance matrix} \\									
\hline									
{} & \multicolumn{2}{c|}{n=50} & \multicolumn{2}{c|}{n=100}\\									
\hline									
& p=500     & p=1000  & p=500   & p=1000 \\									
\midrule									
Glasoo	&	1.597(0.109)	&	0.879(0.013)	&	1.29(0.847)	&	0.428(0.007)	\\
PDSCE	&	0.587(0.13)	&	0.736(0.014)	&	0.094(0.058)	&	0.288(0.01)	\\
CLIME	&	10.5(0.535)	&	11.5(0.233)	&	10.5(0.563)	&	11.5(0.245)\\	
JPEN	&	0.551(0.075)	&	0.691(0.008)	&	0.066(0.042)	&	0.201(0.007)	\\
\bottomrule									
\multicolumn{5}{ |c| }{Toeplitz type covariance matrix} \\									
\hline									
{} & \multicolumn{2}{c|}{n=50} & \multicolumn{2}{c|}{n=100}\\									
\hline									
& p=500     & p=1000  & p=500   & p=1000 \\									
\midrule									
Glasoo	&	2.862(0.475)	&	2.89(0.048)	&	2.028(0.267)	&	2.073(0.078)	\\
PDSCE	&	1.223(0.5)	&	1.238(0.065)	&	0.49(0.269)	&	0.473(0.061)	\\
CLIME	&	4.91(0.22)	&	7.597(0.34)	&	5.27(1.14)	&	8.154(1.168)\\	
JPEN	&	1.151(0.333)	&	2.718(0.032)	&	0.607(0.196)	&	2.569(0.057)	\\
\bottomrule									
\multicolumn{5}{ |c| }{Cov-I type covariance matrix} \\									
\hline									
{} & \multicolumn{2}{c|}{n=50} & \multicolumn{2}{c|}{n=100}\\									
\hline									
& p=500     & p=1000  & p=500   & p=1000 \\									
\midrule									
Glasoo	&	54.0(0.19)	&	190.(5.91)	&	14.7(0.37)	&	49.9(0.08)	\\
PDSCE	&	28.8(0.19)	&	45.8(0.32)	&	16.9(0.04)	&	26.9(0.08)	\\
CLIME	&	59.8(0.82)	&	207.5(3.44)	&	15.4(0.03)	&	53.7(0.69)	\\
JPEN	&	26.3(0.36)	&	7.0(0.07)	&	15.7(0.08)	&	23.5(0.3)	\\
\bottomrule									
\end{tabular}%
\label{tab:addlabel}%
\end{table}%
The Ledoit-Wolf estimate was obtained using code from (http://econ.uzh.ch/faculty/wolf/publications.html\#9). The PDSC estimate was obtained using PDSCE package (http://cran. r-project. org/web/ packages/PDSCE/index.html). The Bickel and Levina's estimator was computed as per the algorithm given in their paper. For inverse covariance matrix performance comparison we include glasso, CLIME (\citet{Cai1:6}) and PDSCE. For each of covariance and inverse covariance matrix estimate, we calculate Average Relative Error (ARE) based on 50 iterations using following formula,
\begin{eqnarray*}
ARE (\Sigma,\hat{\Sigma})= |log(f(S,\hat{\Sigma}))~-~log(f(S,\Sigma_0))|/|(log(f(S,\Sigma_0))|,
\end{eqnarray*}
\begin{figure}[t]
\caption{\textit{Heatmap of zeros identified in covariance matrix out of 50 realizations. White color is 50/50 zeros identified, black color is 0/50 zeros identified.} }
  \begin{center}
 \includegraphics[width=0.85\textwidth]{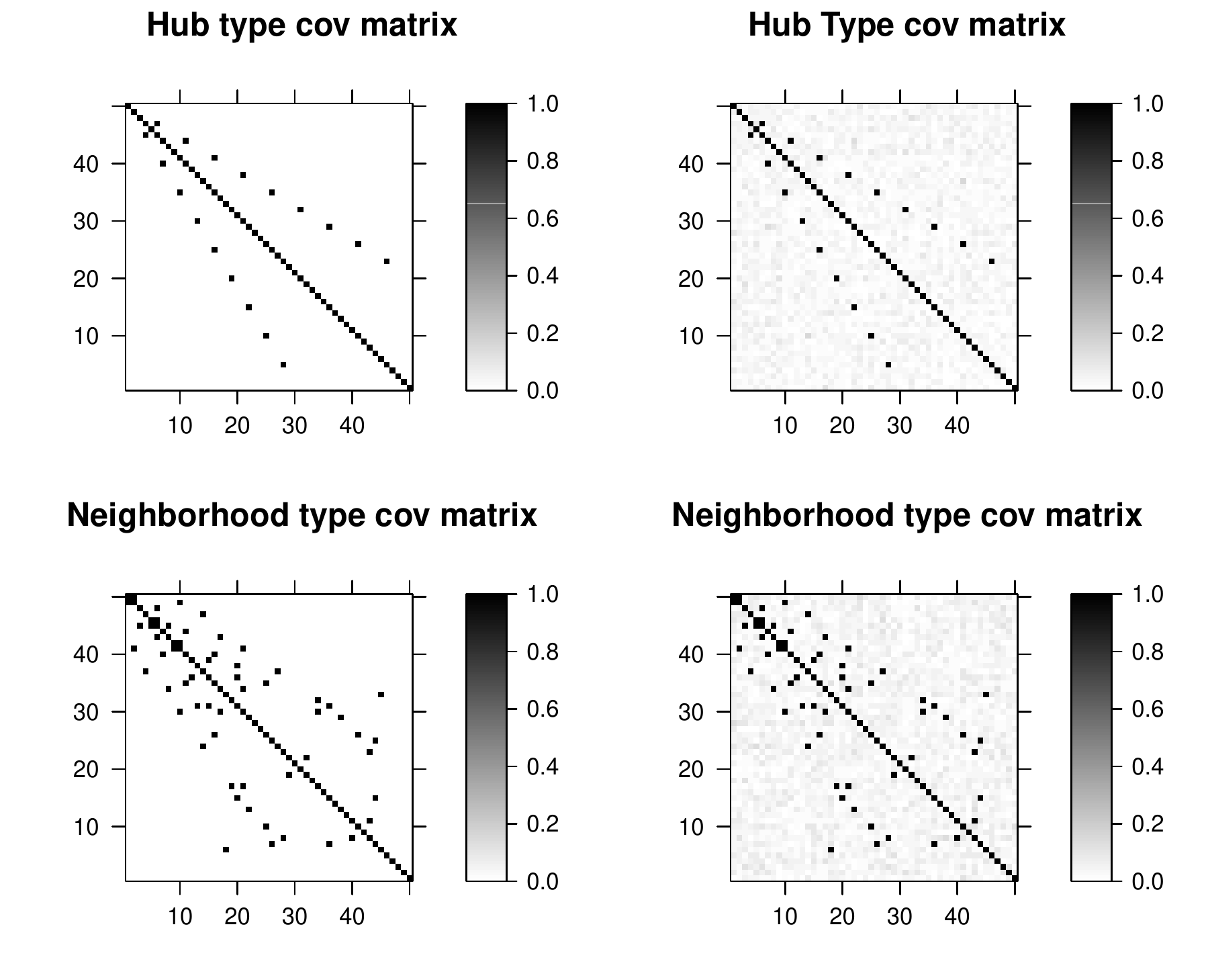}
\end{center}
\end{figure}
where $ f(S,\cdot) $ is multivariate normal density given the sample covariance matrix $S$, $\Sigma_0$ is the true covariance, $\hat{\Sigma}$ is the estimate of $\Sigma_0$ based on one of the methods under consideration. Other choices of performance criteria are Kullback-Leibler used by \citet{Yuan:34} and \citet{Bickel1:4}. The optimal values of tuning parameters were obtained over a grid of values by minimizing $5-$fold cross-validation as explained in $\S4$. The average relative error and their standard deviations (in percentage) for covariance and inverse covariance matrix estimates are given in Table 5.1 and Table 5.2, respectively. The numbers in the bracket are the standard errors of relative error based on the estimates using different methods. Among all the methods JPEN and PDSCE perform similar for most of choices of $n$ and $p$ for all five type of covariance matrices. This is due to the fact that both PDSCE and JPEN use quadratic optimization function with a different penalty function. The behavior of Bickel and levina's estimator is quite good in Toepltiz case where it performs better than the other methods. For this type of covariance matrix, the entries away from the diagonal decay to zero and therefore soft-thresholding estimators like BLThresh perform better in this setting. However for neighorhood and hub type covariance matrix which are not necessarily banded type, Bickel and Levina estimator is not a natural choise as their estimator would fail to recover the underlying sparsity pattern. The performance of Ledoit-Wolf estimator is not very encouraging for Cov-I type matrix. The Ledoit-Wolf estimator shrinks the sample covariance matrix towards identity and hence the eigenvalues estimates are highly shrunk towards one. This is also visible in eigenvalues plot in Figure 5.2 and Figure 5.3. For Cov-I type covariance matrix where most of eigenvalues are close to zero and widely spread, the performance of JPEN estimator is impressive. The eigenplot in Figure 5.3 shows that among all the methods, estimates of eigenvalues of JPEN estimator are most consistent with true eigenvalues. This clearly shows the advantage of JPEN estimator of covariance matrix when the true eigenvalues are dispersed or close to zero. The eigenvalues plot in Figure 5.2 shows that when eigen-spectrum of true covariance matrix are not highly dispersed, the JPEN and PDSCE estimates of eigenavlues are almost the same. This phenomenon is also apparent in Figure 2.1. Also Ledoit-Wolf estimator heavily shrinks the eigenvalues towards the center and thus underestimates the true eigen-spectrum. 
\begin{figure}[t]
  \caption{\textit{Eigenvalues plot for $n=100, p=50$ based on 50 realizations for neighborhood type of covariance matrix}}
  \begin{center}
  \includegraphics[width=.8\textwidth]{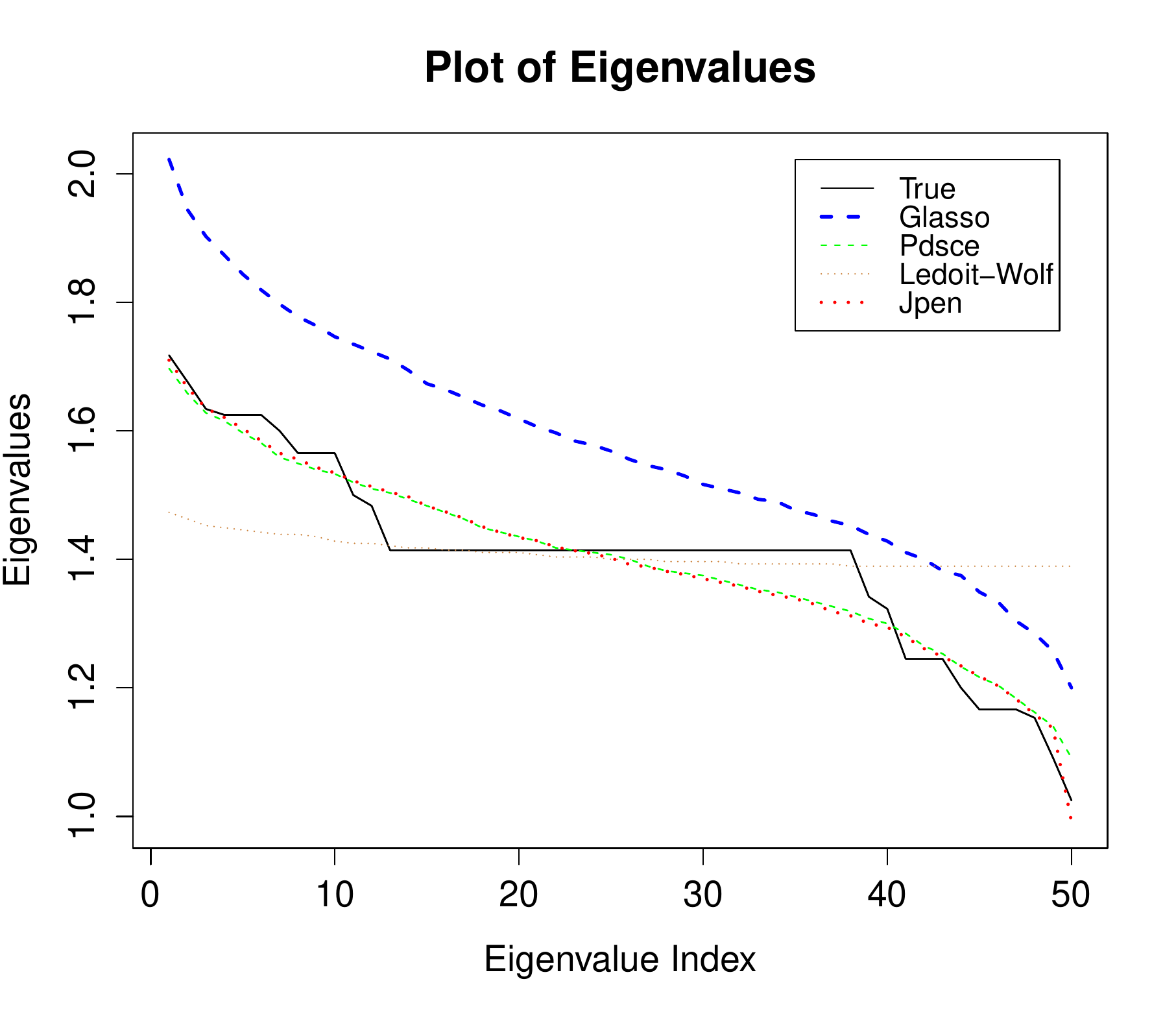}
  \end{center}
\end{figure}
\begin{figure}[t]
  \caption{\textit{Eigenvalues plot for $n=100, p=100$ based on 50 realizations for Cov-I type matrix}}
  \begin{center}
  \includegraphics[width=.8\textwidth]{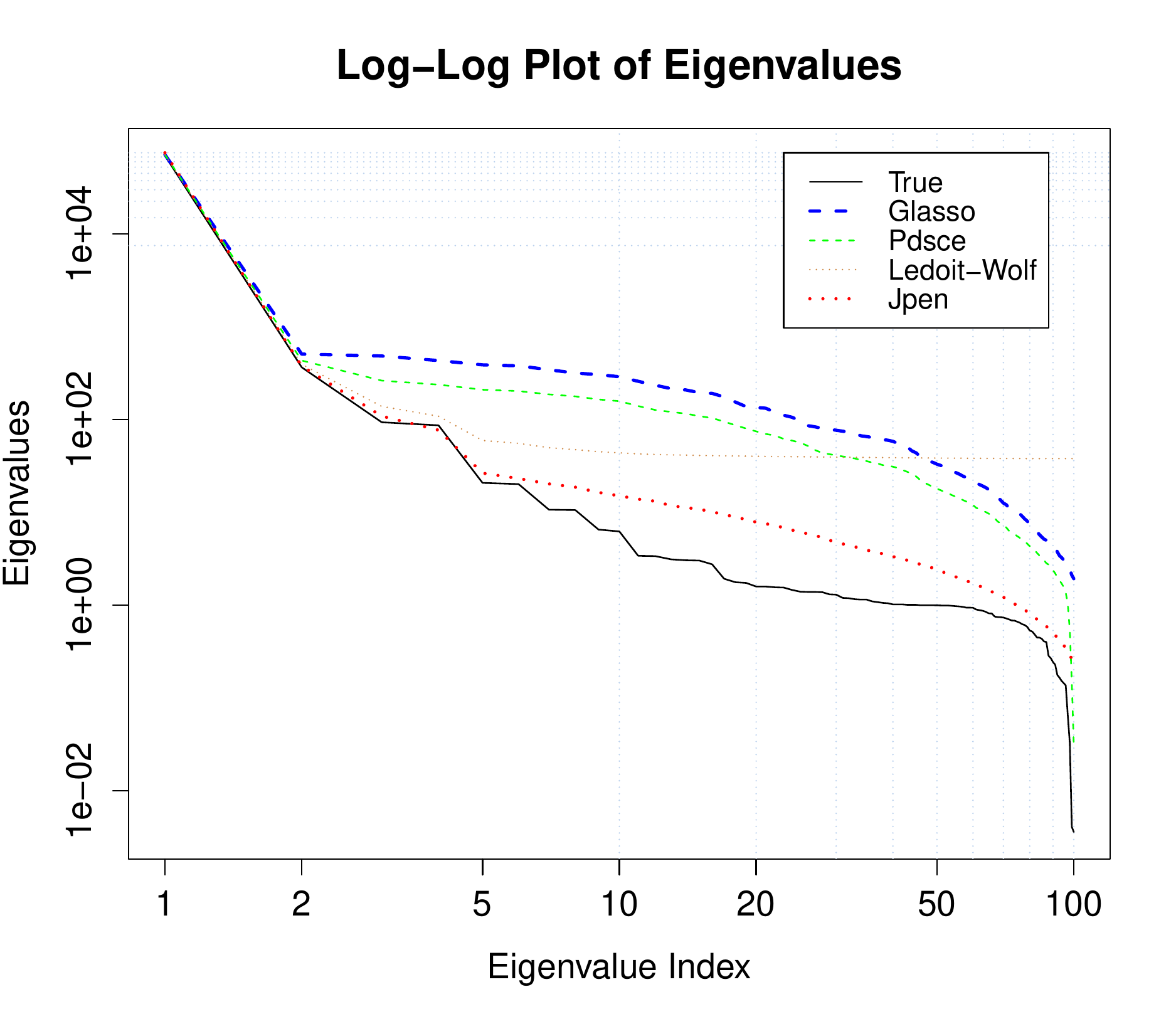}
  \end{center}
\end{figure}

For inverse covariance matrix, we compare glasso, CLIME and PDSCE estimates with proposed JPEN estimator. The JPEN estimator $\hat{\Omega}_{K}$ outperforms other methods for the most of the choices of $n$ and $p$ for all five types of inverse covariance matrices. Additional simulations (not included here) show that for $n \approx p$, all the underlying methods perform similarly and the estimates of their eigenvalues are also well aligned with true values. However in high dimensional setting, for large $p$ and small $n$, their performance is different as seen in simulations of Table 5.1 and Table 5.2. Figure 5.1 shows the recovery of non-zero and zero entries of true covariance matrix based on JPEN estimator $\hat{\Sigma}_{K}$ based on 50 realizations. The estimtor recovers the true zeros for about 90\% of times for Hub and Neighborhood type of covariance matrix. It also reflect the recovery of true structure of non-zero entries and actual pattern among the rows/columns of covariance matrix.
To see the implication of eigenvalues shrinkage penalty as compared to other methods, we plot (Figure 5.2) the eigenvalues of estimated covariance matrix for $n=100$, $p=50$ for neighborhood type of covariance matrix. The JPEN estimates of eigen-spectrum are well aligned with true ones and closest being PDSC estimates of eigenvalues. Figure 5.3 shows the recovery of eigenvalues based on estimates using different methods for Cov-I type covariance matrix. For this particular simulation, the eigenvalues are choosen differently than the one described in (v) of $\S5$. The eigenvalues of true covariance matrix are taken to be very diverse with maximum about $10^6$ and smallest eigenvalue about $10^{-6}$. For Cov-I type of matrix, JPEN estimates of eigenvalues are better than other methods.
\section{Colon Tumor Classification Example}
In this section, we compare performance of JPEN estimator of inverse covariance matrix for tumors classification using Linear Discriminant Analysis (LDA). The gene expression data (\citet{Alan:2} consists of 40 tumorous and 22 non-tumorous adenocarcinoma tissue. After preprocessing, data was reduced to a subset of 2,000 gene expression values with the largest minimal intensity over the 62 tissue samples (source: \textit{http://genomics-pubs.princeton.edu/oncology /affydata/index.html}). In our analysis, we reduced the number of genes by selecting $p$ most \text{sign}ificant genes based on logistic regression. We obtain estimates of inverse covariance matrix for $p=50,100,200$ and then use LDA to classify these tissues as either tumorous or non-tumorous (normal). We classify each test observation x to either class k = 0 or k = 1 using the LDA rule
\begin{eqnarray}
  \delta_k(x)=\argmax_k  \Big \{ x^T\hat{\Omega}\hat{\mu_k}~-~\frac{1}{2} \hat{\mu_k}\hat{\Omega}\hat{\mu_k} ~+~log(\hat{\pi}_k) \Big\},
\end{eqnarray}
where $\hat{\pi}_k$ is the proportion of class $k$ observations in the training data, $\hat{\mu}_k$ is the sample mean for class k on the training data, and $\hat{\Omega}:=\hat{\Sigma}^{-1}$ is an estimator of the inverse of the common covariance matrix on the training data computed. Tuning parameters $\lambda$ and $\gamma$ were chosen using 5-fold cross validation. To create training and test sets, we randomly split the data into a training and test set of sizes 42 and 20 respectively; following the approach used by Wang et al. (2007), the training set has 27 tumor samples and 15 non-tumor samples. We repeat the split at random 100 times and measure the average classification error.
\begin{table}[h]
\caption{\textit{Averages and  standard errors  of classification errors over 100 replications in \%.}}
\begin{center}
\begin{tabular}{l | {l} {l}  {c}r}
\hline
Method			& p=50 & p=100 & p=200  \\
\hline
Logistic Regression	& 21.0(0.84) & 19.31(0.89) & 21.5(0.85) \\
SVM	& 16.70(0.85) & 16.76(0.97) & 18.18(0.96) \\
Naive Bayes	& 13.3(0.75) & 14.33(0.85) & 14.63(0.75) \\
Graphical Lasso	& 10.9(1.3) & 9.4(0.89) & 9.8(0.90) \\
Joint Penalty	& \textbf{9.9(0.98)} & \textbf{8.9(0.93)} & \textbf{8.2(0.81)}
\end{tabular}
\end{center}
\end{table}
Since we do not have separate validation set, we do the 5-fold cross validation on training data. At each split, we divide the training data into 5 subsets (fold) where 4 subsets are used to estimate the covariance matrix and one subset is used to measure the classifier's performance. For each split, this procedure is repeated  5 times by taking one of the 5 subsets as validation data. An optimal combination of $\lambda $ and $\gamma$ is obtained by minimizing the $5$-fold cross validation error.\\

The average classification errors with standard errors over the 100 splits are presented in Table 6.1. Since the sample size is less than the number of genes, we omit the inverse sample covariance matrix as it is not well defined and instead include the naive Bayes' and support vector machine classifiers. Naive Bayes has been shown to perform better than the sample covariance matrix in high-dimensional settings (\citet{Bickel2:39}. Support Vector Machine (SVM) is another popular choice for high dimensional classification tool. Among all the methods covariance matrix based LDA classifiers perform far better than Naive Bayes, SVM and Logistic Regression. For all other classifiers the classification performance deteriorates for increasing $p$. For larger $p$, i.e., when more genes are added to the data set, the classification performance of JPEN estimate based LDA classifier initially improves but it deteriorates for large $p$. For $p=2000$, the classifier based on inverse covariance matrix has accuracy of $30\%$. This is due to the fact that as dimension of covariance matrix increases, the estimator does not remain very informative.
\section{Summary}
We have proposed and analyzed regularized estimation of large covariance and inverse covariance matrix using joint penalty. The proposed JPEN estimators are optimal under spectral norm for underlying classs of sparse and well-conditioned covariance and inverse covariance matrices. We also establish its theoretical consistency in Frobenius norm. One of its biggest advantage is that the optimization carries no computational burden and and the resulting algorithm is very fast and easily scalable to large scale data analysis problems. The extensive simulation shows that the proposed estimators performs well for a number of structured covariance and inverse covariance matrices. Also when the eigenvalues of underlying true covariance matrix are highly dispersed, it outperforms other methods (based on simulation analysis). The JPEN estimator recovers the sparsity pattern of the true covariance matrix and provides a good approximation of the underlying eigen-spectrum and hence we expect that PCA will be one of the most important application of the method. Although the proposed JPEN estimators of covariance and inverse covariance matrix do not require any assumption on the structure of true covariance and inverse covariance matrices respectively, any prior knowledge of structure of true covariance matrix might be helpful to choose a suitable weight matrix and hence improve estimation. \\
\acks{The author would like to express the deep gratitude to Professor Hira L. Koul for his valuable and constructive suggestions during the planning and development of this research work. The author would like to thank Dr. Adam Rothman for his valuable discussion and suggestion. The author would also like to thank the two anonymous referees and the action editor Dr. Jie Peng for insightful reviews that helped to improve the original manuscript substantially.}
\vskip 0.2in
\bibliography{sample}
\newpage
\appendix
\numberwithin{equation}{section}
\section*{Appendix A.}
{\bf Proof of Lemma 3.1}\\
Let 
\begin{equation}
f(R)=\|R-K\|^2+\lambda \|R^-\|_1+\gamma \sum_{i=1}^{p} \{\sigma_i(R)-\bar{\sigma}_R\}^2.
\end{equation}
where $\bar{\sigma}_R$ is the mean of eigenvalues of $R$. Due to the constraint $tr(R)=p$, we have $\bar{\sigma}_R=1$. The third term of (.1) can be written as\\
 $$\sum_{i=1}^{p} \{\sigma_i(R)-\bar{\sigma}_R\}^2=tr(R^2)-2~tr(R)+p$$
We obtain,
\begin{align}
\begin{split}
f(R)&= tr(R^2)-2~tr(RK)+tr(K^2)+\lambda\|R^-\|_1  + \gamma\{tr(R^2) -2~tr(R)+p\} \\
&= tr(R^2(1+\gamma))-2~tr(K+\gamma~I)+tr(K^2)+\lambda\|R^-\|_1 + p \\
&= (1+\gamma)\|R-(K+\gamma~I)/(1+\gamma)\|^2+tr(K^2)+\lambda \|R^-\|_1 + p
\end{split}
\end{align}
This is quadratic in $R$ with a $\ell_1$ penalty to the off-diagonal entries of $R$, therefore a convex function in $R$. \\\\
 {\bf Proof of Lemma 3.2}
The solution to (.2) satisfies:
\begin{eqnarray}
2(R-(K+\gamma I))(1+\gamma)^{-1}+ \lambda \frac{\partial{\|R^-\|_1}}{\partial{R}}=0
\end{eqnarray}
where $\frac{\partial{\|R^-\|_1}}{\partial{R}}$ is given by:
\[\frac{\partial{\|R^-\|_1}}{\partial{R}} = \left\{ \begin{array}{lr} 1 & : if~~~R_{ij}>0 \\ -1 &: if ~~~ R_{ij} <0 \\  \tau \in (-1,1)& : if~~~ R_{ij}=0 \end{array} \right. \]
Note that $\|R^-\|_1$ has same value irrespective of \text{sign} of $R$, therefore the right hand side of (.2) is minimum if :\\
\begin{eqnarray*}
\text{sign}(R)=\text{sign}(K+\gamma I)=\text{sign}(K)
\end{eqnarray*}
$\forall \epsilon >0$, using (.3), $\sigma_{min}\{ (K+\gamma I)-\frac{\lambda}{2} \text{sign}(K) \} >\epsilon $ gives a $(\lambda,\gamma) \in \hat{\mathscr{S}}^{K}_1$ and such a choice of $(\lambda,\gamma)$ gaurantees the estimator to be positive definite. \\
{\bf Remark:} Intuitively, a larger $\gamma$ shrinks the eigenvalues towards center which is 1, a larger $\gamma$ would result in positive definite estimator, whereas a larger $\lambda$ results in sparse estimate. A combination of $(\lambda, \gamma)$ results in a sparse and well-conditioned estimator. In particular case, when $K$ is diagonal matrix, the $\lambda<2*\gamma$. \\ \\
{\bf Proof of Theorem 3.1}
Define the function $Q(.)$ as following:
	$$ Q(R)=f(R)-f(R_0)$$
where  $R_0$ is the true correlation matrix and $R$ is any other correlation matrix. Let $R=UDU^T$ be eigenvalue decomposition of $R$, $D$ is diagonal matrix of eigenvalues and $U$ is matrix of eigenvectors.
We have,
\begin{align}
\begin{split}
Q(R)& =  \|R-K\|_F^2+\lambda \|R^-\|_1+\gamma~ tr(D^2-2~D+p) \\
& ~~~-\|R_0-K\|_F^2 - \lambda\|R_0^-\|_1 -\gamma~ tr(D_0^2-2~D_0+p)
\end{split}
\end{align}
$R_0=U_0D_0U_0^T$ is eigenvalue decomposition of $R_0$. Let $ \Theta_n(M):=\{\Delta : \Delta=\Delta^T , ~ \|\Delta\|_2=Mr_n,~ 0< M < \infty ~\}$. The estimate $\hat{R}$ minimizes the $Q(R)$ or equivalently $ \hat{\Delta}=\hat{R}-R_0$ minimizes the $G(\Delta)=Q(R_0+\Delta)$. Note that $G(\Delta)$ is convex and if $\hat{\Delta}$ be its solution, then  we have $G(\hat{\Delta}) \le G(0) = 0$. Therefore if we can show that $G(\Delta)$ is non-negative for $\Delta \in \Theta_n(M)$, this will imply that the $\hat{\Delta}$ lies within sphere of radius $Mr_n$. We require $r_n=o\Big (\sqrt{(p+s)~log~p/n}\Big )$. 
\begin{eqnarray*}
\|R-K\|_F^2-\|R_0-K\|_F^2 & =& tr(R^TR -2 R^TK+K^TK)-tr(R^T_0R_0-2R_0S+K^TK) \\
& = & tr(R^TR-R^T_0R_0)-2~tr((R-R_0)^TK) \\
& = & tr((R_0+\Delta)^T(R_0+\Delta)-R^T_0R_0)-2~tr(\Delta^TK)\\
& = & tr(\Delta^T\Delta)-2~tr(\Delta^T(K-R_0))
\end{eqnarray*}
Next, we bound term involving $K$ in above expression, we have
\begin{eqnarray*}
|tr(\Delta^T(R_0-K))| & \le & \sum_{i\neq j} |\Delta_{ij}({R_0}_{ij}-K_{ij})| \\
& \le & \max_{i\neq j}(|{R_0}_{ij}-K_{ij}|) \|\Delta^-\|_1 \\
& \le&  C_0 (1+\tau) \sqrt{\frac{log~p}{n}} \|\Delta^-\|_1 \le C_1 \sqrt{\frac{log~p}{n}} \|\Delta^-\|_1
\end{eqnarray*}
holds with high probability by a result (Lemma 1) from \citet{Ravi1:24} on the tail inequality for sample covariance matrix of sub-gaussian random vectors and where $C_1=C_0 (1+\tau), C_0 >0$.
Next we obtain upper bound on the terms involving $\gamma$ in (.4). we have,
\begin{eqnarray*}
    \lefteqn{tr(D^2-2D) - tr(D_0^2-2D_0) } \\
    & = & tr\{R^2-R_0^2\} - 2~tr\{R-R_0)\} = tr(R_0+\Delta)^2-tr(R^2_0) \\
			& = & 2~tr(R_0 \Delta) +tr(\Delta^T\Delta) \le 2~\sqrt{s}\|\Delta\|_F+ \|\Delta\|^2_F.
    \end{eqnarray*}
using Cauchy-Schwarz inequality. To bound the term $ \lambda(\|R^-\|_1-\|R_0^-\|_1) =\lambda(\|\Delta^-+R_0^-\|_1-\|R_0^-\|_1) $, let $E$ be index set as defined in Assumption A.2 of Theorem 3.2. Then using the triangle inequality, we obtain,
\begin{eqnarray*}
\lambda(\|\Delta^-+R_0^-\|_1-\|R_0^-\|_1) & = & \lambda(\|\Delta_E^-+R_{0}^-\|_1+\|\Delta_{\bar{E}}^-\|_1-\|R_0\|_1) \\
& \geq & \lambda(\|R_{0}^-\|_1-\|\Delta_E^-\|_1+\|\Delta_{\bar{E}}^-\|_1 -\|R_0^-\|_1) \\
& \geq & \lambda(\|\Delta_{\bar{E}}^-\|_1-\|\Delta_{E}^-\|_1)
\end{eqnarray*}
Let $ \lambda =(C_1/\epsilon)\sqrt{log~p/n}$, $\gamma=(C_1/\epsilon_1)\sqrt{log~p/n}, $ where $(\lambda, \gamma) \in \hat{\mathscr{S}}^{K}_1$, we obtain,
\begin{eqnarray*}
G(\Delta)& \geq & tr(\Delta^T\Delta)(1+\gamma)-2~C_1 \Big \{ \sqrt{\frac{log~p} {n}}(\|\Delta^{-}\|_1)+ \frac{1}{\epsilon_1}\sqrt{\frac{s~log~p} {n}} \|\Delta\|_F \Big \} \\
& &  +\frac{C_1}{\epsilon} \sqrt{\frac{log~p}{n}} \big ( \|\Delta^-_{\bar{E}}\|_1  - \Delta^-_{E}\|_1 \big )\\
& \geq & \|\Delta\|_F^2(1+\gamma) -2C_1\sqrt{\frac{log~p}{n}} \big ( \|\Delta^-_{\bar{E}}\|_1 + \|\Delta^-_{{E}}\|_1 \big ) \\
& & \frac{C_1}{\epsilon} \sqrt{\frac{log~p}{n}} \big ( \|\Delta^-_{\bar{E}}\|_1- \Delta^-_{E}\|_1 \big )  - \frac{2C_1}{\epsilon_1}\sqrt{\frac{s~log~p} {n}} \|\Delta\|_F.
\end{eqnarray*}
Also because $\|\Delta^-_{E}\|_1=\sum_{(i,j) \in E, i \neq j } \Delta_{ij} \leq \sqrt{s} \|\Delta^-\|_F$,
\begin{eqnarray*}
 -2C_1\sqrt{\frac{log~p}{n}} \|\Delta^-_{\bar{E}}\|_1 +\frac{C_1}{\epsilon} \sqrt{\frac{log~p}{n}} {\|\Delta^-_{\bar{E}}\|}_1 & \geq & \sqrt{\frac{log~p}{n}} \|\Delta^-_{\bar{E}}\|_1 \big ( -2 C_1 + \frac{C_1}{\epsilon} \big )  \\
& \geq & 0
\end{eqnarray*}
for sufficiently small $\epsilon$. Therefore,
\begin{eqnarray*}
G(\Delta) & \geq  & \|\Delta\|_F^2 \big ( 1+\frac{C_1}{\epsilon_1}\sqrt{\frac{log~p}{n}} \big ) -C_1\sqrt{\frac{s~log~p}{n}} \|\Delta^+\|_F \{1+1/\epsilon+2/\epsilon_1\}\\
& \geq & \|\Delta\|_F^2\Big[1 +\frac{C_1}{\epsilon_1}\sqrt{\frac{log~p}{n}}-\frac{C_1}{M}\{1+1/\epsilon+2/\epsilon_1\} \Big] \\
 & \geq & 0,
\end{eqnarray*}
for all sufficiently large $n$ and $M$. Which proves the first part of theorem. To prove the operator norm consistency, we have,
\begin{eqnarray*}
\|\hat{\Sigma}_K-\Sigma_0\| & = & \|\hat{W}\hat{R}\hat{W}-WK W\| \\
& \le & \|\hat{W}-W\| \|\hat{R}-K\| \|\hat{W}-W\| \\
&  & + \|\hat{W}-W\| (\|\hat{R}\|\|W\| +\|\hat{W}\| \|K\|) +  \|\hat{R}-K\| \|\hat{W}\| \|W\|.
\end{eqnarray*}
using sub-multiplicative norm property $\|AB\| \leq \|A\|\|B\|$. Since $\|K\|=O(1)$ and $\|\hat{R}-K\|_F=O(\sqrt{\frac{s~log~p}{n}})$ these together implies that $\|\hat{R}\|=O(1)$ . Also,
\begin{eqnarray*}
\|{\hat{W}}^2-W^2\| & = & \max_{\|x\|_2=1} \sum_{i=1}^p |({\hat{w}_i}^2-w_i^2)| x_i^2 \leq \max_{1 \leq i \leq p} |({\hat{w}_i}^2-w_i^2)| \sum_{i=1}^p x_i^2 \\
& = & \max_{1 \leq i \leq p} |({\hat{w}_i}^2-w_i^2)| = O \big ( \sqrt{\frac{log~p}{n}} \big ).
\end{eqnarray*}
holds with high probability by using a result (Lemma 1) from \citet{Ravi1:24}. Next we shall shows that $\|\hat{W}-W\|\asymp \|\hat{W}^2-W^2\|$, (where A$\asymp$B means A=$O_P(B)$ and B=$O_P(A)$). We have,
\begin{eqnarray*}
\|\hat{W}-W\| & = & \max_{\|x\|_2=1} \sum_{i=1}^p |({\hat{w}_i}-w_i)| x_i^2 = \max_{\|x\|_2=1} \sum_{i=1}^p | \big (\frac{{{\hat{w}_i}}^2-w_i^2}{\hat{w}_i+w_i} \big ) | x_i^2 \\
 & \asymp & \sum_{i=1}^p |({\hat{w}_i}^2-w_i^2)| x_i^2 = C_3 \|\hat{W}^2-W^2\|.
\end{eqnarray*}
where we have used the fact that the true standard deviations are well above zero, i.e., $\exists~ 0 < C_3 < \infty $ such that $1/C_3 \leq w^{-1}_i \leq C_3 ~\forall ~i=1,2,\cdots, p$, and sample standard deviation are all positive, i.e, $\hat{w}_i > 0 ~\forall ~i=1,2,\cdots,p.$ Now since $\|\hat{W}^2-W^2\| \asymp \|\hat{W}-W\|$, this follows that $\|\hat{W}\|=O(1)$ and we have $\|\hat{\Sigma}_K-\Sigma_0\|^2=O \big(\frac{s~log~p}{n} +\frac{log~p}{n} \big )$. This completes the proof. \\ \\
\textbf{Proof of Theorem 3.2}
Let
\begin{eqnarray*}
f(\Sigma)=||\Sigma-S||^2_F +  \lambda \|\Sigma^-\|_1  + \gamma\sum_{i=1}^{p} \{\sigma_i(\Sigma)-\bar{\sigma}_{\Sigma}\}^2,
\end{eqnarray*}
Similar to the proof of theroem (3.1), define the function $Q_1(.)$ as following:
$$ Q_1(\Sigma)=f(\Sigma)-f(\Sigma_0)$$
where  $\Sigma_0$ is the true covariance matrix and $\Sigma$ is any other covariance matrix. Let $\Sigma=UDU^T$ be eigenvalue decomposition of $\Sigma$, $D$ is diagonal matrix of eigenvalues and $U$ is matrix of eigenvectors.
We have,
\begin{align}
\begin{split}
Q_1(\Sigma)& =  \|\Sigma-S\|_F^2+\lambda \|\Sigma^-\|_1+\gamma~ tr(D^2)-(tr(D))^2/p \\
& ~~~-\|\Sigma_0-S\|_F^2 - \lambda\|\Sigma_0^-\|_1 -\gamma~ tr(D_0^2)-(tr(D_0))^2/p
\end{split}
\end{align}
where $A=diag(a_1,a_2,\cdots,a_p)$ and $\Sigma_0=U_0D_0U_0^T$ is eigenvalue decomposition of $\Sigma_0$. Write $\Delta=\Sigma-\Sigma_0$, and let $ \Theta_n(M):=\{\Delta : \Delta=\Delta^T , ~ \|\Delta\|_2=Mr_n,~ 0< M < \infty ~\}$. The estimate $\hat{\Sigma}$ minimizes the $Q(\Sigma)$ or equivalently $ \hat{\Delta}=\hat{\Sigma}-\Sigma_0$ minimizes the $G(\Delta)=Q(\Sigma_0+\Delta)$. Note that $G(\Delta)$ is convex and if $\hat{\Delta}$ be its solution, then  we have $G(\hat{\Delta}) \le G(0) = 0$. Therefore if we can show that $G(\Delta)$ is non-negative for $\Delta \in \Theta_n(M)$, this will imply that the $\hat{\Delta}$ lies within sphere of radius $Mr_n$. We require $\sqrt{(p+s)~log~p}=o \Big (\sqrt{n} \Big ) $.
\begin{eqnarray*}
\|\Sigma-S\|_F^2-\|\Sigma_0-S\|_F^2 & =& tr(\Sigma ^T \Sigma -2 \Sigma^TS+S^TS)-tr(\Sigma_0^T\Sigma_0-2\Sigma_0S+S^TS) \\
& = & tr(\Sigma^T\Sigma-\Sigma_0^T\Sigma_0)-2~tr((\Sigma-\Sigma_0)S) \\
& = & tr((\Sigma_0+\Delta)^T(\Sigma_0+\Delta)-\Sigma_0^T\Sigma_0)-2~tr(\Delta^TS)\\
& = & tr(\Delta^T\Delta)-2~tr(\Delta^T(S-\Sigma_0))
\end{eqnarray*}
Next, we bound term involving $S$ in above expression, we have
\begin{eqnarray*}
|tr(\Delta(\Sigma_0-S))| & \le & \sum_{i\neq j} |\Delta_{ij}({\Sigma_0}_{ij}-S_{ij})| + \sum_{i=1} |\Delta_{ii}({\Sigma_0}_{ii}-S_{ii})| \\
& \le & \max_{i\neq j}(|{\Sigma_0}_{ij}-S_{ij}|) \|\Delta^-\|_1 + \sqrt{p} \max_{i=1} (|{\Sigma_0}_{ii}-S_{ii}|) \sqrt{\sum_{i=1} \Delta^2_{ii}} \\
& \le&  C_0 (1+\tau) \max_{i}(\Sigma_{0ii}) \Big \{ \sqrt{\frac{log~p}{n}} \|\Delta^-\|_1+ \sqrt{\frac{p~log~p}{n}} \|\Delta^+\|_2 \Big \} \\
& \le & C_1 \Big \{ \sqrt{\frac{log~p}{n}} \|\Delta^-\|_1+ \sqrt{\frac{p~log~p}{n}} \|\Delta^+\|_2 \Big \}
\end{eqnarray*}
holds with high probability by a result (Lemma 1) from \citet{Ravi1:24} where $C_1=C_0 (1+\tau) \max_{i}(\Sigma_{0ii}), C_0 >0$ and $\Delta^+$ is matrix $\Delta$ with all off-diagonal elements set to zero. Next we obtain upper bound on the terms involving $\gamma$ in (3.7). we have,
 \begin{eqnarray*}
    \lefteqn{tr(D^2)-(tr(D))^2/p - tr(D_0^2)-(tr(D))^2/p } \\
    & = & tr(\Sigma^2)-tr(\Sigma_0^2) - (tr(\Sigma))^2/p + (tr(\Sigma_0))^2/p
    \end{eqnarray*}
\begin{eqnarray*}
\textbf{(i)}~~~\lefteqn{tr(\Sigma^2)-\Sigma_0^2))} \\
& \leq &  tr(\Sigma_0+\Delta)^2 -tr(\Sigma_0)^2 \\
& = & tr(\Delta)^2+ 2~tr(\Delta^2\Sigma_0) \leq tr(\Delta)^2 + C_1 \sqrt{s}\|\Delta\|_F
\end{eqnarray*}
\begin{eqnarray*}
\textbf{(ii)}~~~\lefteqn{tr((\Sigma))^2 -(tr(\Sigma_0))^2}\\
& = & (tr(\Sigma_0+\Delta))^2 - (tr(\Sigma_0))^2 \\
& \leq & (tr(\Delta))^2 +2~tr(\Sigma_0)~tr(\Delta) \leq p~\| \Delta\|^2_F +2~\bar{k} p\sqrt{p} \|\Delta^+\|_F.
\end{eqnarray*}
Therefore the $\gamma$ term can be bounded by $2\|\Delta\|^2_F+(C_1\sqrt{s}+2\sqrt{p}\bar{k})\|\Delta\|_F$. We bound the term invloving $\lambda$ as in similar to the proof of Theorem 3.1. For $\lambda \asymp \gamma \asymp \sqrt{\frac{log~p}{n}}$, the proof follows very simialr to Therem 3.1.\\ \\
{\bf Proof of Theorem 3.3.}
To bound the cross product term involving $\Delta$ and $\hat{R}_K^{-1}$, we have,
\begin{eqnarray*}
|tr((R_0^{-1}-\hat{R}_K^{-1})\Delta)| & = & |tr(R_0^{-1}(\hat{R}_K-R_0)\hat{R}_K^{-1}{\Delta})| \\
& \leq & \sigma_1(R_0^{-1})|tr((\hat{R}_K-R_0)\hat{R}_K^{-1} \Delta)| \\
& \leq & \bar{k}\sigma_1(\hat{R}_K^{-1})|tr((\hat{R}_K-R_0) \Delta)| \\
& \leq & \bar{k} \bar{k}_1|tr((\hat{R}_K-R_0)\Delta)|.
\end{eqnarray*}
where $\sigma_{min}(\hat{R}_K) \geq (1/\bar{k}_1) >0 $, is a positive lower bound on the eigenvalues of JPEN estimate $\hat{R}_K$ of correlation matrix $R_0$. Such a constant exist by Lemma 3.2. Rest of the proof closely follows as that of Theorem 3.1. \\ \\
{\bf Proof of Theorem 3.4.}
We bound the term $tr((\hat{\Omega}_{S}-\Omega_0)\Delta)$ similar to that in proof of Theorem 3.3. Rest of the proof closely follows to that Theorem 3.2.

\end{document}